\documentclass[sigconf,screen]{acmart}

\usepackage{cmll}

\usepackage{mathtools} 
\usepackage{stmaryrd}
\usepackage{mathpartir}
\usepackage{tikz,tikz-cd}
\usetikzlibrary{arrows,arrows.meta,positioning}
\usepackage{url}
\usepackage{enumerate,xspace}

\usepackage[utf8]{inputenc}
\usepackage[english]{babel}
\usepackage[draft]{optional}
\usepackage{comment}
\usepackage[colorinlistoftodos,prependcaption,textsize=tiny]{todonotes}
\usepackage{xcolor}
\usepackage[all]{xy}
\usepackage{braket}
\usetikzlibrary{fit}

\usepackage{graphicx}
\graphicspath{  }

\usepackage{cmll}

\newcommand{\Ob}{\text{Ob}}

\newcommand{\cat}[1]{\mathbb{#1}}

\DeclareMathOperator{\fix}{\mathbf{fix}}

\newcommand{\eval}{\text{eval}} 

\newcommand{\catname}[1]{\mathbf{#1}}
\newcommand{\Rel}{\catname{Rel}}

\newcommand{\Pow}{\catname{PW}}

\newcommand{\Tang}{\mathbf{T}} 

\newcommand{\lis}[1]{\langle #1 \rangle}
\DeclareMathOperator{\Diff}{\mathbf{D}}
\DeclareMathOperator{\Trace}{\mathbf{Tr}} 

\newcommand{\eg}{\emph{e.g.}}

\newcommand{\true}{\mathbf{t}}
\newcommand{\false}{\mathbf{f}}

\newcommand{\One}{\mathbf{1}} 
\newcommand{\Zero}{\mathbf{0}} 
\newcommand{\NN}{\mathbb{N}}

\newcommand{\QQ}{\mathbb{Q}}
\newcommand{\RR}{\mathbb{R}}

\newcommand{\N}{\mathcal{N}}
\newcommand{\Y}{\mathcal{Y}}
\newcommand{\M}{\mathcal{M}}
\newcommand{\Z}{\mathcal{Z}}
\newcommand{\J}{\mathrm{J}}
\newcommand{\Rp}{\mathbb{R}_+}

\newcommand{\ints}[1]{\underline{#1}}
\newcommand{\Sring}{\mathbb{S}} 

\makeatletter
\def\slashedarrowfill@#1#2#3#4#5{%
	$\m@th\thickmuskip0mu\medmuskip\thickmuskip\thinmuskip\thickmuskip
	\relax#5#1\mkern-7mu%
	\cleaders\hbox{$#5\mkern-2mu#2\mkern-2mu$}\hfill
	\mathclap{#3}\mathclap{#2}%
	\cleaders\hbox{$#5\mkern-2mu#2\mkern-2mu$}\hfill
	\mkern-7mu#4$%
}
\def\rightslashedarrowfill@{%
	\slashedarrowfill@\relbar\relbar\mapstochar\rightarrow}
\newcommand\xslashedrightarrow[2][]{%
	\ext@arrow 0055{\rightslashedarrowfill@}{#1}{#2}}
\makeatother
\def\proft{\xslashedrightarrow{}}

\AtBeginDocument{%
  }

\setcopyright{acmlicensed}

\copyrightyear{2024}
\acmYear{2024}
\setcopyright{acmlicensed}\acmConference[LICS '24]{39th Annual ACM/IEEE
	Symposium on Logic in Computer Science}{July 8--11, 2024}{Tallinn, Estonia}
\acmBooktitle{39th Annual ACM/IEEE Symposium on Logic in Computer Science
	(LICS '24), July 8--11, 2024, Tallinn, Estonia}
\acmDOI{10.1145/3661814.3662108}
\acmISBN{979-8-4007-0660-8/24/07}

\begin{document}

\title{Combining fixpoint and differentiation theory}

\author{Zeinab Galal}
\affiliation{%
  \institution{University of Bologna}
  \country{Italy}
}
\affiliation{%
	\institution{Inria Sophia Antipolis}
	\country{France}
}
\author{Jean-Simon Pacaud Lemay}
\affiliation{%
  \institution{Macquarie University}
  \country{Australia}
\email{larst@affiliation.org}
}


\begin{abstract}
Interactions between derivatives and fixpoints have many important applications in both computer science and mathematics. In this paper, we provide a categorical framework to combine fixpoints with derivatives by studying Cartesian differential categories with a fixpoint operator. We introduce an additional axiom relating the derivative of a fixpoint with the fixpoint of the derivative. We show how the standard examples of Cartesian differential categories where we can compute fixpoints provide canonical models of this notion. We also consider when the fixpoint operator is a Conway operator, or when the underlying category is closed. As an application, we show how this framework is a suitable setting to formalize the Newton-Raphson optimization for fast approximation of fixpoints and extend it to higher order languages.
\end{abstract}




\settopmatter{printfolios=true}
\maketitle

\section{Introduction}

Fixpoint theory has many applications for establishing the existence and uniqueness of solutions for systems of differential and integral equations. Dually, using the information provided by derivatives to approximate more efficiently fixpoints is a key tool in numerical analysis, data-flow analysis, automatic differentiation, enumerative combinatorics, formal language theory, etc. A well-known example is the Newton-Raphson method for fast approximation of roots or fixpoints of differentiable functions. Starting from an initial point, the main idea is to compute the next iteration by using the tangent line. It is a powerful technique as the convergence rate is quadratic: at each iteration, the number of accurately computed digits doubles. Many other standard optimization methods are based on generalizations or variants of the Newton-Raphson method. 

From the viewpoint of programming languages, fixpoints are a key feature as they allow for the presence of loops and defining programs via recursion or iteration. Incorporating recursive features into differential programming languages has seen increasing interest in recent years \cite{abadi2019simple,mazza2021automatic,di2023language} and there is already a rich literature on implicit automatic differentiation to perform automatic differentiation efficiently on functions defined implicitly by fixpoint equations \cite{gilbert1992automatic, bell2008algorithmic, blondel2022efficient}.

While the categorical semantics of differentiation and fixpoints are now well-established lines of research, one aspect that has not yet been studied in full detail is the interaction between them from a semantical viewpoint. The objective of this paper is to develop a denotational framework to combine fixpoints with derivatives by axiomatizing the notion of Cartesian differential categories with a fixpoint operator. This axiomatization provides a guideline to introducing differentials to $\lambda$-calculi with fixpoint operators such as PCF \cite{plotkin1977lcf} and the $\lambda$-$\mathrm{Y}$-calculus \cite{statman2002lambdaY}. 

Cartesian differential categories (Definition \ref{cartdiffdef}) were introduced by Blute, Cockett, and Seely in \cite{blute2009cartesian}, and provide the categorical foundations of multivariable differential calculus. Cartesian \emph{closed} differential categories (Section \ref{sec:closed}) provide the categorical semantics of the differential $\lambda$-calculus, introduced by Ehrhard and Regnier~\cite{ehrhard2003differential, bucciarelli2010categorical, manzonetto2012categorical}. Cartesian (closed) differential categories have been quite successful in formalizing various notions related to differentiation, and have also found numerous applications in both mathematics and computer science. Moreover, Cartesian differential categories belong to a larger story of categorical foundations of differentiation, which include differential categories~\cite{blute2006differential}, coherent differential categories~\cite{thomas2023cartesian}, reverse differential categories~\cite{cockett2019reverse}, differential restriction categories~\cite{cockett2011differential}, and tangent categories~\cite{cockett2014differential}. 

In this paper, we show how to combine the differential operator in a Cartesian differential category with various categorical notions of fixpoint operators \cite{bloom1993iteration, Hasegawa97recursionfrom, simpson2000complete, hasegawa2004uniformity, joyal1996traced} such as parametrized fixpoint operators (Definition \ref{ParamFixpointCartesian}) and Conway operators (Definition \ref{def: Conway}). As such, to this end, we introduce \emph{Cartesian differential fixpoint categories} (Definition \ref{def:fixdiff}), which are Cartesian differential categories equipped with a parametrized fixpoint operator which satisfies what we call the \emph{differential-fixpoint rule} (\ref{ax:fixdiff}) describing the derivative of the parametrized fixpoint. 

To help understand the differential-fixpoint rule, let us provide some intuition using ordinary calculus. Consider the smooth function $f  : \RR \times \RR \to \RR$ defined as $f(x,y) := 1 - x^2 y$. Then the function $g : \RR \to \RR$ defined as $g(x):= \frac{1}{1+ x^2}$ is a solution of the fixpoint or implicit equation
\begin{equation}
	\label{eqn:fixexample}
	g(x) = f(x, g(x)).
\end{equation}
If we differentiate both sides, using the chain rule on the right-hand side, we obtain:
\begin{equation}
	\label{eqn:fixderexample}
	g'(x) = \frac{\partial f}{\partial x} (x, g(x)) + \frac{\partial f}{\partial y} (x, g(x))  \cdot g'(x)
\end{equation}
Provided that $\frac{\partial f}{\partial y} (x, g(x)) \neq 1$, we may rewrite (\ref{eqn:fixderexample}) as: 
\[
g'(x) = \left(1 - \frac{\partial f}{\partial y} (x, g(x))\right)^{-1} \frac{\partial f}{\partial x} (x, g(x)) = \frac{-2x}{(1+x^2)^2}
\]
relating the derivative of $f$ with the derivative of its fixpoint. This is the basis of many iteration schemes to compute or approximate the derivative of a fixpoint such as the Newton-Raphson method. If we combine the equations (\ref{eqn:fixexample}) and (\ref{eqn:fixderexample}), we obtain that the function $\Tang(g): \RR \times \RR \to \RR \times \RR$ defined as  $\Tang(g)(x,a)= (g(x), g'(x)\cdot a)$ is a fixpoint for the function $\Tang(f):  \RR \times \RR \times \RR \times \RR \to \RR \times \RR$ defined as: 
\begin{equation}
	\label{eq:tangentexample}
	\Tang(f)(x, y, a, b)= (f(x,y), \frac{\partial f}{\partial x} (x, y) \cdot a + \frac{\partial f}{\partial y} (x, y)  \cdot b ).
\end{equation}
This $\Tang$ operation is called the \emph{tangent bundle functor} (Definition~\ref{def:tan-func}), and can be defined in any Cartesian differential category (it corresponds to the tangent bundle in algebraic geometry and to the Jacobian vector product in forward differentiation). Therefore, computing the derivative of the fixpoint of $f$ is equivalent to computing the fixpoint of $\Tang(f)$, and this is the property we axiomatize for Cartesian differential categories.

\textbf{Outline:} Section \ref{sec:Cartdiff} is a background section where we review Cartesian differential categories. In Section \ref{sec:fixdif}, we study the compatibility relation between parametrized fixpoint operators and differential combinators, and introduce Cartesian differential fixpoint categories. In Section \ref{sec:traced}, we consider the case of Conway operators, and prove an equivalence between various axioms relating derivatives and fixpoints. In Section \ref{subsec:fix-linear}, we study the relation between fixpoints and linearity. In Section \ref{sec:closed}, we consider the Cartesian closed setting. In Section \ref{sec:examplesfixdif}, we provide many examples of Cartesian differential fixpoint categories. In Section \ref{sec:Newton}, we provide an application by extending the Newton-Raphson optimization method to our setting. We conclude in Section \ref{sec:conclusion} with a discussion about future work.

\textbf{Related works:} In \cite[Theorem 5.29]{ehrhard2023coherent}, Ehrhard proves a compatibility relation in a cpo-enriched setting between the least fixpoint operator (in the coKleisli category) and the tangent bundle functor of a Scott coherent differential category. In \cite{sprunger2019differentiable}, Sprunger and Katsumata construct Cartesian differential categories with delayed trace operators, which are related to trace operators but no longer assume the fixpoint axiom. 

\textbf{Conventions:} We assume the reader is familiar with basic notions of category theory such as categories, functors, products, etc. An arbitrary category will be denoted by $\cat{C}$, with arbitrary objects denoted by capital letters $A, B, X, Y$, etc. and maps by minuscule letters $f, g, h$, etc. Hom-sets are denoted as $\cat{C}(A,B)$, arbitrary maps as $f: A \to B$, identity maps as $1_A: A \to A$, and we use the standard notation $\circ$ and convention for composition (unlike in other Cartesian differential papers which use diagrammatic order). 

\section{Cartesian Differential Categories}\label{sec:Cartdiff}

In this background section, we review Cartesian differential categories. For a more in-depth introduction, we refer the reader to~\cite{blute2009cartesian, lemay2021exponential, lemay2023cartesian, garner2020cartesian, cockett2020linearizing}. 

The underlying structure of a Cartesian differential category is that of a \emph{Cartesian left additive category}, which can be described as a category with finite products which is \emph{skew}-enriched over the category of commutative monoids \cite{garner2020cartesian}. More explicitly, a \textbf{left additive category} \cite[Definition 1.1.1]{blute2009cartesian} is a category $\cat{C}$ such that each hom-set $\cat{C}(A,B)$ is a commutative monoid (written additively), so we have zero maps $0$ and can take the sum of maps $f+g$, however only pre-composition preserves the additive structure, that is, 
\[
(f+g) \circ x = f \circ x + g \circ x \quad \text{and} \quad 0 \circ x = 0.
\]Maps $f$ which do preserve the additive structure, that is, 
\[
f \circ (x + y) = f \circ x + f \circ y \quad \text{and} \quad f \circ 0 =0,
\] are called \textbf{additive} maps \cite[Definition 1.1.1]{blute2009cartesian}. Then a Cartesian left additive category is a left additive category with finite products which is compatible with the additive structure. For a category $\cat{C}$ with finite products, we denote the product by $\times$, the projection maps as $\pi_j: A_1 \times \hdots \times A_n \to A_j$, the pairing operation by $\langle -, \hdots, - \rangle$, the canonical diagonal map by $\Delta_A = \langle 1_A, 1_A \rangle : A \to A \times A$ and the terminal object by $\top$. Then a \textbf{Cartesian left additive category} \cite[Definition 1.2.1]{blute2009cartesian} is a left additive category $\cat{C}$ which has finite products and such that all the projection maps are additive. 

A \emph{Cartesian differential category} is a Cartesian left additive category equipped with a \emph{differential combinator}, which is an operator which sends maps to their derivative. In this paper, it is important to note that we follow the now more widely used convention used for Cartesian differential categories which flips the convention used in \cite{blute2009cartesian}, so that the linear argument of the derivative is in the second argument rather than in the first. 

\begin{definition}\label{cartdiffdef}
	A \emph{Cartesian differential category} \cite[Definition 2.1.1]{blute2009cartesian} is a Cartesian left additive category $\cat{C}$ equipped with a \textbf{differential combinator} $\Diff$, which is a family of functions indexed by pairs of objects in $\cat{C}$:
	\begin{align*}
		\Diff : \cat{C}(A , B) \to \cat{C}(A \times A, B) && \infer{f: A \to B}{\Diff[f]: A \times A \to B}
	\end{align*}
	such that the seven following axioms hold: 
	\begin{description}
		\item[{\bf [CD.1]}] $\Diff[0]=0$ and for all $f, g : A \to B$, \[\Diff[f + g] = \Diff[f] +\Diff[g]\]
		\item[{\bf [CD.2]}] for all $f : A \to B$ and $x,y,z: C \to A$,
		\[\begin{aligned}
		\Diff[f] \circ \langle x, 0 \rangle &= 0 \quad \text{and}\\
		\Diff[f] \circ \langle x,y+ z \rangle &= \left( \Diff[f] \circ \langle x, y \rangle \right) +  \left( \Diff[f] \circ \langle x, z \rangle \right)
		\end{aligned}\]
		\item[{\bf [CD.3]}] $\Diff[1_A]=\pi_2$ and for $j \in \{1,2\}$, \[\Diff[\pi_j] = \pi_{j} \circ \pi_2 : (A_1 \times A_2) \times (A_1 \times A_2)  \to A_j\]
		\item[{\bf [CD.4]}] for all $f_i : A \to B_i$,
		\[\Diff[\left\langle f_1, \hdots, f_n \right \rangle] = \left \langle  \Diff[f_1], \hdots, \Diff[f_n] \right \rangle\]
		\item[{\bf [CD.5]}] for all $f: A \to B$, $g: B\to C$,
		\[\Diff[g \circ f] = \Diff[g] \circ \langle f \circ \pi_1, \Diff[f] \rangle\]
		\item[{\bf [CD.6]}] for all $f : A \to B$ and $x,y,z: C \to A$,
		\[\Diff\left[\Diff[f] \right] \circ \left \langle x, y, 0, z \right \rangle=  \Diff[f] \circ \langle x, z \rangle\]
		\item[{\bf [CD.7]}] for all $f : A \to B$ and $x,y,z: C \to A$,
		\[\Diff\left[\Diff[f] \right] \circ \left \langle x, y, z, 0 \right \rangle = \Diff\left[\Diff[f] \right] \circ \left \langle x, z, y, 0 \right \rangle\]
	\end{description}
	For a map $f: A \to B$, the map $\Diff[f]: A \times A \to B$ is called the \textbf{derivative} of $f$. 
\end{definition}

The axioms of a differential combinator are analogues of the basic properties of the total derivative from multivariable differential calculus. The axioms say that: (1) the derivative of a sum is the sum of the derivatives, (2) derivatives are additive in their second argument, (3) the derivative of identity maps and projections are projections, (4) the derivative of a pairing is the pairing of the derivatives, (5) the chain rule for the derivative of a composition, (6) the derivative is \emph{linear} in its second argument, and lastly (7) is the symmetry of the partial derivatives. The term linear refers to the canonical notion of linearity in a Cartesian differential category, which we discuss in Section \ref{subsec:fix-linear}. There is a sound and complete term calculus for Cartesian differential categories \cite[Section 4]{blute2009cartesian}, which is useful for intuition and proofs. So we write:
\[ \Diff[f](a,b) := \frac{\mathsf{d}f(x)}{\mathsf{d}x}(a) \cdot b\]
In particular, the chain rule axiom is expressed as: 
\[ \frac{\mathsf{d}g\left(f(x) \right)}{\mathsf{d}x}(a) \cdot b = \frac{\mathsf{d}g(y)}{\mathsf{d}y}(f(a)) \cdot \left( \frac{\mathsf{d}f(x)}{\mathsf{d}x}(a) \cdot b \right) \]

Arguably, the canonical example of a Cartesian differential category is the category of real smooth functions, whose differential combinator is given by the standard differentiation of smooth functions \cite[Example 2.10]{cockett2020linearizing}. In Section \ref{sec:closed}, we will discuss Cartesian \emph{closed} differential categories, which are particularly important since they provide the categorical semantics of the differential $\lambda$-calculus~\cite{bucciarelli2010categorical} -- though closed structure does not play a crucial role for the story of this paper. In Section \ref{sec:examplesfixdif}, we provide other examples of Cartesian differential categories. For more examples of Cartesian differential categories, we refer the reader to see \cite{cockett2020linearizing,garner2020cartesian}. 

In a Cartesian differential category, the differential combinator induces a functor called the tangent bundle functor. This functor plays an important role in the story of this paper since, as we will see in Section \ref{sec:fixdif}, it fits very naturally with fixpoint operators. The name comes from the fact that every Cartesian differential category is a \textbf{tangent category} \cite[Proposition 4.7]{cockett2014differential}, and so this functor corresponds to the classical tangent bundle functor for Euclidean spaces. 

\begin{definition}\label{def:tan-func} For a Cartesian differential category $\cat{C}$, the \textbf{tangent bundle functor} \cite[Proposition 1]{lemay2021exponential} is the endofunctor $\Tang: \cat{C} \to \cat{C}$ defined on objects as $\Tang (A) = A \times A$ and on maps as $\Tang (f) = \left \langle f \circ \pi_1, \mathsf{D}[f] \right \rangle$, which in the term calculus is expressed as: 
	\[\mathsf{T}(f)(a,b) = \left(f(a), \frac{\mathsf{d}f(x)}{\mathsf{d}x}(a) \cdot b \right)\] 
\end{definition}

Properties that the tangent bundle functor satisfies can be found in \cite[Lemma 2]{lemay2021exponential}. It is worthwhile to point out that the tangent bundle preserves composition if and only if the chain rule of the differential combinator holds. In fact, the chain rule can be reformulated in terms of the tangent functor bundle as: $\Diff[g\circ f] = \Diff[g] \circ \Tang(f)$.

\section{Combining Fixpoints and Derivatives}\label{sec:fixdif}

In this section, we combine differential combinators and parametri-zed fixpoint operators by studying their compatibility in an arbitrary Cartesian differential category, and then introduce the notion of a Cartesian differential fixpoint category. For a more in-depth introduction to categorical interpretations of fixpoint operators, we refer the reader to \cite{bloom1993iteration, Hasegawa97recursionfrom, simpson2000complete, hasegawa2004uniformity}.  

Let us begin by explaining why for a Cartesian differential category, one must consider a parametrized fixpoint operator rather than simply a basic fixpoint operator. Recall that in a category with a terminal object, a \textbf{point} of an object $X$ is a map from the terminal object to $X$, so $p: \top \to X$. Then a \textbf{fixpoint operator} \cite[Definition 2.1]{simpson2000complete} is an operator $\fix$ which for every endomorphsim $f: X \to X$ associates a point $\fix(f): \top \to X$ such 
\[
f \circ \fix(f) = \fix(f)
\]
meaning that $\fix(f)$ is a fixpoint of $f$. Unfortunately, this kind of fixpoint operator is in a certain sense incompatible with differential combinators. This is because, in a Cartesian differential category, the derivative of a point $p: \top \to X$ is always zero, $\Diff[p] = 0$. Therefore: 
\begin{lemma}\label{lem:diff-fix-0} In a Cartesian differential category with a fixpoint operator, for every map $f: X \to X$, $\Diff[\fix(f)] = 0$. 
\end{lemma}

For a Cartesian differential category, one must instead consider \emph{parametrized} fixpoint operators, which axiomatize the notion of fixpoints for maps in context. For a map of type $A \times X \to X$, the parameter $A$ is viewed as representing the context of the term, then taking the parametrized fixpoint gives a map of type $A \to X$. Parametrized fixpoint operators are axiomatized by two axioms: the fixpoint axiom and by naturality in the context argument. 

\begin{definition}\label{ParamFixpointCartesian} For a category $\cat{C}$ with finite products, a \textbf{parame-trized fixpoint operator} \cite[Definition 2.2]{simpson2000complete} is a family of functions $\fix$ indexed by pairs of objects in $\cat{C}$, 
	\begin{align*}
		\fix^X_A: \cat{C}(A \times X, X) \to \cat{C}(A, X) && \infer{f: A \times X \to X}{  \fix^X_A(f): A  \to X}
	\end{align*}
	such that the following axioms hold: \\
	\noindent 1. \textbf{Parametrized fixpoint axiom:} for all maps ${f :A \times X\to X}$:
	\[f \circ  \lis{1_A,\fix^X_A (f)} = \fix^X_A (f)\]
	\noindent 2. \textbf{Naturality axiom:} for all maps $g : A\to B$ and $f :B \times  X \to X$:
	\[\fix^X_B (f) \circ g = \fix^X_A(f \circ (g \times 1_X))\]
	For a map $f: A \times X \to X$, the map $ \fix^X_A(f): A  \to X$ is called the \textbf{parametrized fixpoint} of $f$. 
\end{definition}

Since we will be working with term calculus notation for Cartesian differential categories, we shall use the following term calculus notation to write the parametrized fixpoint operator: 
\[ \fix^X_A(f)(a) = \mu x.f(a,x) \]
where the variable $x$ is bound. The notation $\mu$ here is for an arbitrary fixpoint operator and not necessarily the least fixpoint operator.
For example, the parametrized fixpoint axiom, which says that $\fix^X_A(f)$ is a fixpoint of $f$ in context $A$, is expressed as: 
\begin{equation}
	\label{eqn:fixaxiom}
	f(a, \mu x. f(a,x)) =\mu x. f(a,x).
\end{equation}
A well-known example of a category with a parametrized fixpoint operator is the category of Scott domains and Scott continuous functions \cite[Example 7.1.2]{Hasegawa97recursionfrom}, whose parametrized fixpoint operator is given by the standard Kleene iteration formula. In Section \ref{sec:closed}, we will discuss parametrized fixpoint operators in Cartesian \emph{closed} categories. In Section \ref{sec:examplesfixdif}, we provide other examples of parametrized fixpoint operators. For more examples of parametrized fixpoint operators, we refer the reader to see \cite{Hasegawa97recursionfrom,hasegawa2004uniformity}.

So how should a differential combinator and parametrized fixpoint operator interact? In particular, what should the derivative of a parametrized fixpoint be? Consider the following computation using the parametrized fixpoint axiom and the chain rule:
\begin{equation} 
	\begin{split}
	&\;\;\frac{\mathsf{d} \mu x. f(u,x)}{\mathsf{d}u}(a) \cdot b = \frac{\mathsf{d} f \left( u, \mu x. f(u,x) \right)}{\mathsf{d}u}(a) \cdot b\\
	=&\; \frac{\mathsf{d}f(u,v)}{\mathsf{d}(u,v)} \left( a, \mu x. f(a,x) \right) \cdot ( b, \frac{\mathsf{d} \mu x. f(u,x)}{\mathsf{d}u}(a) \cdot b ) 
	\label{eqn:derfix}
\end{split}
\end{equation}
From $(\ref{eqn:fixaxiom})$ and $(\ref{eqn:derfix})$, we see that the pair 
\[
(\mu x. f(a,x), \frac{\mathsf{d} \mu x. f(u,x)}{\mathsf{d}u}(a) \cdot b)
\]is a parametrized fixpoint of $\Tang(f)(a,x,b,y)$ in the $x$ and $y$ variables, that is: 
\begin{gather*} 
	\Tang(f)(a, \mu x. f(a,x), b, \frac{\mathsf{d} \mu x. f(u,x)}{\mathsf{d}u}(a) \cdot b ) \\
	= (\mu x. f(a,x), \frac{\mathsf{d} \mu x. f(u,x)}{\mathsf{d}u}(a) \cdot b)
\end{gather*}
Thus, our compatibility relation between a parametrized fixpoint operator and differential combinator is asking that the derivative of the parametrized fixpoint be equal to the second component of the nested fixpoint: 
\[\mu(x,y).(f(a,x), \frac{\mathsf{d}f(u,v)}{\mathsf{d}(u,v)}(a,x) \cdot (b,y)).\] 
In other words, the derivative of the parametrized fixpoint is the second component of the parametrized fixpoint of the tangent bundle, up to a rearranging of variables. To make certain that this indeed makes sense, let us check that the types work correctly. So starting with a map ${f: A \times X \to X}$, taking its parametrized fixpoint gives us $\fix^X_A(f): A \to X$, and so its derivative is of type $\Diff[\fix^X_A(f)]: A \times A \to X $. On the other hand, applying the tangent bundle functor gives us $\Tang(f): A \times X \times A \times X \to X \times X$. For this to be of the correct type to apply the parametrized fixpoint operator, we must swap the middle two terms. So let 
\[c = \langle \pi_1, \pi_3, \pi_2, \pi_4 \rangle : A \times B \times C \times D \to A \times C \times B \times D\] be the canonical isomorphism which swaps the middle two terms. Then pre-composing by $c$ gives us $\Tang(f) \circ c: A \times A \times X \times X \to X \times X$. So we may take its parametrized fixpoint to get a map of type $\fix_{A \times A}^{X \times X}\left(\Tang(f)\circ  c \right): A \times A \to X \times X$. Lastly, we post-compose by the second projection to finally get 
\[\pi_2 \circ \fix_{A \times A}^{X \times X}\left(\Tang(f)\circ  c \right): A \times A \to X.\] When these two maps are equal (Figure \ref{fig:fixderaxiom}), we call this the differential-fixpoint rule and use it as our axiom for the definition of a Cartesian differential fixpoint category.

\begin{figure*}[h]
	\fbox{\begin{minipage}{.65\textwidth}\begin{center}
				\vspace{-.5em}
				\[  \begin{array}[c]{c} \inferrule *
					{\inferrule * {  A \times X \xrightarrow{f} X }{ A \xrightarrow{\fix(f)} X }}
					{A \times A \xrightarrow{\Diff[\fix(f)]} X}
				\end{array} \qquad=\qquad \begin{array}[c]{c}
					\inferrule *
					{\inferrule *
						{\inferrule * {  A \times X \xrightarrow{f} X }{ A \times X \times A \times X \xrightarrow{\Tang(f)} X \times X}}
						{ A \times A \times X \times X \xrightarrow{c} A \times X \times A \times X \xrightarrow{\Tang(f)} X \times X}}
					{A \times A \xrightarrow{\fix(\Tang(f)\circ c))} X \times X \xrightarrow{\pi_2} X} \end{array} \]
				\vspace{-.5em}
	\end{center}\end{minipage}}
	\caption{\normalfont Differential-fixpoint rule\label{fig:fixderaxiom}}
\end{figure*}

\begin{definition}\label{def:fixdiff} A parametrized fixpoint operator for a Cartesian differential category satisfies the \textbf{differential-fixpoint rule} if for every map $f: A \times X \to X$, the following equality holds: 
	\begin{align}
		\label{ax:fixdiff}
		\Diff [ \fix^X_A(f)] = \pi_2 \circ \fix_{\Tang A}^{\Tang X}(\Tang(f)\circ c)
	\end{align}
	which in the term calculus is expressed as follows: 
	\[
	\frac{\mathsf{d} \mu x. f(u,x)}{\mathsf{d}u}(a) \cdot b = 
	\pi_2 ( \mu(x,y).(f(a,x), \frac{\mathsf{d}f(u,v)}{\mathsf{d}(u,v)}(a,x) \cdot (b,y) ) )
	\]
	A \textbf{Cartesian differential fixpoint category} is a Cartesian differential category with a parametrized fixpoint operator which satisfies the differential-fixpoint rule. 
\end{definition}

We can also ask how the tangent bundle functor and the parame-trized fixpoint operator interact. From the differential-fixpoint rule, by definition we have that: 
\[ \Tang ( \fix^X_A(f)) = \langle \fix^X_A(f) \circ \pi_1, \pi_2 \circ \fix_{\Tang A}^{\Tang X}\left(\Tang(f)\circ c\right)  \rangle \]
However, since the first component of $\fix_{\Tang A}^{\Tang X}\left(\Tang(f)\circ c\right)$ is not necessarily $\fix^X_A(f) \circ \pi_1$, we may not have that $\Tang ( \fix^X_A(f))$ is equal to $\fix_{\Tang A}^{\Tang X}\left(\Tang(f)\circ c\right)$. When these two are equal, we call this the tangent-fixpoint rule. 

\begin{definition}\label{def:tangentfix} A parametrized fixpoint operator for a Cartesian differential category satisfies the \textbf{tangent-fixpoint rule} if for every map $f: A \times X \to X$, the following equality holds: 
	\begin{align}
		\label{ax:fixtang}
		\Tang ( \fix^X_A(f) ) =  \fix_{\Tang A}^{\Tang X}(\Tang(f)\circ c)
	\end{align}
	which in the term calculus is expressed as follows: 
\begin{gather*} 
	(\mu x. f(a,x), \frac{\mathsf{d} \mu x. f(u,x)}{\mathsf{d}u}(a) \cdot b )\\ = 
	\mu(x,y).(f(a,x), \frac{\mathsf{d}f(u,v)}{\mathsf{d}(u,v)}(a,x) \cdot (b,y) ) 
\end{gather*} 
	It means that the Cartesian functor $(\Tang,c)$ is a morphism of categories with parametrized fixpoint operator (Definition 1.3 in \cite{bloom1993iteration} in the coCartesian setting).
\end{definition}

It is straightforward to see that the tangent-fixpoint rule implies the differential-fixpoint rule: 

\begin{lemma}\label{lem:tanfix-difffix} A parametrized fixpoint operator which satisfies the tangent-fixpoint rule also satisfies the differential-fixpoint rule. 
\end{lemma}  
\begin{proof} By definition, note that $\pi_2 \circ \Tang(f) = \Diff(f)$. Therefore, post-composing both sides of (\ref{ax:fixtang}) by the second projection $\pi_2$ results precisely in (\ref{ax:fixdiff}).    
\end{proof}

Therefore, a Cartesian differential category with a parametrized fixpoint operator which satisfies the tangent-fixpoint rule is a Cartesian differential fixpoint category. While we do not have an example of a Cartesian differential category with fixpoint operator where the differential-fixpoint rule holds but not the tangent-fixpoint rule, there is no reason to assume that the converse is necessarily true. The reason for choosing the differential-fixpoint rule over the tangent-fixpoint rule as an axiom for Cartesian differential fixpoint categories comes from the fact that it is the differential combinator which is central in the definition of Cartesian differential categories rather than the tangent bundle functor. When asking for compatibility regarding a structure or property on a Cartesian differential category, it is more natural to ask what the derivatives of the necessary structural or property maps are. The tangent-fixpoint rule is more natural when considering a Cartesian differential category instead as a tangent category \cite{cockett2014differential}, and will be the central axiom for the notion of a tangent category with a parametrized fixpoint operator (which is future work, see Section \ref{sec:conclusion}). That being said, we show that for \emph{Conway} fixpoint operators these two rules are equivalent.

\section{Conway Operators and Derivatives}\label{sec:traced}

Fixpoint operators are closely related to the notion of \emph{trace} operators. In particular, to give a trace operator for products corresponds to providing a special kind of parametrized fixpoint operator called a \emph{Conway} operator. In this section, we study the compatibility between differential combinators and Conway operators (and trace operators), and introduce the notion of a traced Cartesian differential category. For a more in-depth introduction to Conway operators and trace operators, we refer the reader to \cite{joyal1996traced,Hasegawa97recursionfrom,hasegawa2004uniformity}.

\begin{definition} For a category $\cat{C}$ with finite products, a \textbf{trace operator} \cite[Section 2]{joyal1996traced} is a family of functions $\Trace$ indexed by triples of objects in $\cat{C}$, 
	\begin{align*}
		\Trace^X_{A,B}: \cat{C}(A \times X, B \times X) \to \cat{C}(A, B) && \infer{f: A \times X \to B \times X}{ \Trace^X_{A,B}(f): A  \to X}
	\end{align*}
	which satisfy the axioms found in \cite[Definition 2.1]{Hasegawa97recursionfrom}. 
\end{definition}

As mentioned, for a category with finite products, giving a trace operator is equivalent to giving a special kind of parametrized fixpoint operator called a Conway operator.  We provide here the axiomatization of a Conway operator as found in \cite{Hasegawa97recursionfrom}, in terms of dinaturality and Beki\v{c}'s axiom. Equivalent axiomatizations can be found in \cite{hasegawa2004uniformity,simpson2000complete}. 

\begin{definition}\label{def: Conway} For a category $\cat{C}$ with finite products, a \textbf{Conway operator} \cite[Theorem 3.1]{Hasegawa97recursionfrom} is a parametrized fixpoint operator $\fix$ which also satisfies: \\
	\noindent 1. \textbf{Dinaturality axiom:} for all maps $f : A \times X \to Y$ and $g : Y \to X$
	\begin{gather*}
		\fix^X_A(g \circ f) = g \circ \fix^Y_A (f \circ (1_A \times g))
	\end{gather*}
	which in the term calculus is expressed as:
	\[ \mu x.g \left( f(x,a) \right) = g\left( \mu x. f(x, g(a)) \right) \]
	\noindent 2. \textbf{Beki\v{c}'s axiom:} for all maps $f: A \times X \times Y \to X$ and $g: A \times X \times Y \to Y$ the following equality holds: 
	\begin{gather*}
		\fix^{X \times Y}_A(\langle f, g \rangle) =\\
		\langle l, \fix^Y_A(g \circ (1_A \times l \times 1_Y)\circ (\Delta_A \times 1_Y)) \rangle
	\end{gather*}
where $l := \fix^X_A(f \circ ((1_{A\times X} \times \fix^Y_{A \times X}(g)) \circ \Delta_{A \times X})):A \to X$. In the term calculus, it is expressed as:
	\begin{gather*} \mu (x,y) \left( f(a,x,y), g(a,x,y) \right) =\\
		\left( \mu x. f\!\left(a,x, \mu y. g(a,x,y) \right), \mu y. g \!\left( a,\mu x. f\!\left(a,x, \mu y. g(a,x,y) \right)\!, y \right) \right)
	\end{gather*}
\end{definition}

While the right-hand side of Beki\v{c}'s axiom's expression seems more complex, it is extremely useful in practice. The idea is that Beki\v{c}'s axiom can be understood as a form of Gaussian elimination: solving a system of nested fixpoint equations with multiple variables is reduced to solving fixpoint equations with one variable at a time. As such, Beki\v{c}'s axiom will be key in the proofs of this section. Most known fixpoint operators are Conway but it is possible to construct fixpoint operators that satisfy all of the axioms but one of the axioms in Definitions \ref{ParamFixpointCartesian} and \ref{def: Conway} \cite{esik1988independence}.

For a category with finite products, there is a bijective correspondence between trace operators and Conway operators \cite[Theorem 3.1]{Hasegawa97recursionfrom}. To go from a Conway operator to a trace operator, first recall that by the universal property of the product, a map of type $f: A \times X \to B \times X$ is a tuple of maps $f = \langle f_1, f_2 \rangle$ where $f_1 = \pi_1 \circ f: A \times X \to B$ and $f_2= \pi_2 \circ f:  A \times X \to X$. Then starting with a Conway operator $\fix$, define the trace operator $\Trace$ on a map $f = \langle f_1, f_2 \rangle: A \times X \to B \times X$ as the following composite: 
\[ \Trace^X_{A,B}( f )(a) =  f_1( \lis{ 1_A, \fix^X_A(f_2)} ) \]
which in the term calculus is expressed as: 
\[\Trace^X_{A,B}\left( \langle f_1, f_2 \rangle \right)(a) = f_1 ( a, \mu x. f_2(a,x) ). \]
On the other hand, starting from a trace operator $\Trace$, define the Conway operator $\fix$ on a map $f: A \times X \to X$ as follows: 
\[\fix^X_{A}(f) = \mathsf{Tr}^X_{A,X}\left(\langle f, f \rangle \right)  \]
As such, a category with finite products and a Conway operator is a \textbf{traced (Cartesian) monoidal category} \cite{joyal1996traced}. So we define: 

\begin{definition}
	A \textbf{traced Cartesian differential category} is a Cartesian differential fixpoint category whose parametrized fixpoint operator is a Conway operator. 
\end{definition}

We wish to show that for a Conway operator, the differential-fixpoint rule is equivalent to the tangent-fixpoint rule. To do so, it will be useful to consider yet another equivalent rule that even more precisely describes the derivative of the parametrized fixpoint. Notice that when we applied the chain rule to $\frac{\mathsf{d} \mu x. f(u,x)}{\mathsf{d}u}(a) \cdot b$ above, we in fact showed that it was a parametrized fixed point for \[
\frac{\mathsf{d}f(u,v)}{\mathsf{d}(u,v)} \left( a, \mu x. f(a,x) \right) \cdot \left( b, - \right).
\] We may then ask that it instead be \emph{the} parametrized fixpoint. 

\begin{definition} A parametrized fixpoint operator for a Cartesian differential category satisfies the \textbf{strong differential-fixpoint rule} if for every map $f: A \times X \to X$, the following equality holds: 
	\begin{align}\label{eq:strongrule}
		\Diff[ \fix^X_A(f) ] = \fix^X_{A \times A} ( \Diff[f] \circ \langle \pi_1, \fix^X_A(f) \circ \pi_1, \pi_2, \pi_3  \rangle ) 
	\end{align}
	which in the term calculus is expressed as follows: 
	\begin{align*}
		\frac{\mathsf{d} \mu x. f(u,x)}{\mathsf{d}u}(a) \cdot b = \mu y. \frac{\mathsf{d}f(u,v)}{\mathsf{d}(u,v)} \left( a, \mu x. f(a,x) \right) \cdot \left( b, y \right)    
	\end{align*}
\end{definition}

We should explain why we chose the differential-fixpoint rule as the axiom for Cartesian differential fixpoint categories instead of the strong differential-fixpoint rule. While the strong differential-fixpoint rule only involves the differential combinator, notice that $\fix^X_A(f)$ appears on both sides of (\ref{eq:strongrule}). Thus, while both the differential combinator and the tangent bundle functor are used in the differential-fixpoint rule, it is slightly more natural. With all that said, we will show that for a Conway operator, thanks to Beki\v{c}'s axiom, all three of these rules are equivalent. 

\begin{proposition}\label{prop:tanfix-difffix} For a Conway operator $\fix$ of a Cartesian differential category, the following are equivalent: 
	\begin{enumerate}[{\em (i)}]
		\item $\fix$ satisfies the differential-fixpoint rule;
		\item $\fix$ satisfies the tangent-fixpoint rule;
		\item $\fix$ satisfies the strong differential-fixpoint rule. 
	\end{enumerate}
\end{proposition}
\begin{proof} The key to this proof is Beki\v{c}'s axiom. Before we try to prove the equivalence of these three rules, let us first compute a useful identity using the term calculus. Assuming we have a Conway operator (and not assuming any of the three rules), consider any map $f: A \times X \to X$. To simplify using Beki\v{c}'s axiom, let $z = (a,b)$, and define 
	\[h(z,x,y) = f(a,x) \quad \text{and} \quad k(z,x,y) = \frac{\mathsf{d}f(u,v)}{\mathsf{d}(u,v)}(a,x) \cdot (b,y).\]Then using Beki\v{c}'s axiom, we compute that: 
	\begin{gather*} 
		 \mu(x,y).\left( h(z,x,y), k(z,x,y) \right) \\
		= ( \mu x. h(z,x, \mu y. k(z,x,y) ), \mu y. k ( z,\mu x. h(z,x, \mu y. k(z,x,y) ), y ) ) 
	\end{gather*} 
	So we have the equality:
	\begin{equation}\label{eq:special}\begin{gathered} 
			\mu(x,y).(f(a,x), \frac{\mathsf{d}f(u,v)}{\mathsf{d}(u,v)}(a,x) \cdot (b,y) ) \\
			= ( \mu x.f(a,x), \mu y. \frac{\mathsf{d}f(u,v)}{\mathsf{d}(u,v)} ( a, \mu x. f(a,x) ) \cdot \left( b, y \right) ) 
		\end{gathered}
	\end{equation}
	Using this identity, we can now prove the desired equivalence. We will prove $(ii) \Rightarrow (i) \Rightarrow (iii) \Rightarrow (ii)$. Now $(ii) \Rightarrow (i)$ is just Lemma \ref{lem:tanfix-difffix}. For $(i) \Rightarrow (iii)$, using (\ref{eq:special}), we can easily expand out the differential-fixpoint rule as: 
	\begin{gather*}
		\frac{\mathsf{d} \mu x. f(u,x)}{\mathsf{d}u}(a) \cdot b  \\
		=  \pi_2 ( \mu(x,y).(f(a,x), \frac{\mathsf{d}f(u,v)}{\mathsf{d}(u,v)}(a,x) \cdot (b,y) ) ) \\
		= \pi_2( \mu x.f(a,x), \mu y. \frac{\mathsf{d}f(u,v)}{\mathsf{d}(u,v)} ( a, \mu x. f(a,x) ) \cdot \left( b, y \right) ) \\
		= \mu y. \frac{\mathsf{d}f(u,v)}{\mathsf{d}(u,v)} ( a, \mu x. f(a,x) ) \cdot ( b, y )
	\end{gather*}
	So the strong differential-fixpoint rule holds. For $(iii) \Rightarrow (ii)$, applying the strong differential-fixpoint rule to (\ref{eq:special}) gives us that: 
	\begin{gather*}
		\mu(x,y).(f(a,x), \frac{\mathsf{d}f(u,v)}{\mathsf{d}(u,v)}(a,x) \cdot (b,y) ) \\
		= ( \mu x.f(a,x), \mu y. \frac{\mathsf{d}f(u,v)}{\mathsf{d}(u,v)} ( a, \mu x. f(a,x) ) \cdot ( b, y ) ) \\
		=  \mu(x,y).(f(a,x),\frac{\mathsf{d} \mu x. f(u,x)}{\mathsf{d}u}(a) \cdot b) 
	\end{gather*}
	So the tangent-fixpoint rule holds. So we conclude that for a Conway operator, the three rules are equivalent. 
\end{proof}

Therefore, for a traced Cartesian differential category, its Conway operator satisfies the differential-fixpoint rule, the tangent-fixpoint rule, and the strong differential-fixpoint rule. Since a Conway operator is equivalent to a trace operator, there is yet another equivalent rule involving the trace operator. Unfortunately, there does not appear to be a nice formula for describing the derivative of the trace. Instead, the trace operator behaves quite nicely with the tangent bundle functor. Now for a map $f: A \times X \to B \times X$, we have that \[\Tang(\Trace^X_{A,B}\left( f \right)  ): A \times A \to B \times B.\]On the other hand, first applying the tangent functor gives $\Tang\left( f \right): A \times X \times A \times X \to B \times X \times B \times X$. Before we can take its trace, we must post- and pre-compose by $c$ to get $c \circ \Tang( f ) \circ c: A \times A \times X \times X \to B \times B \times X \times X$. Then taking the trace we finally get a map of type
\[
\Trace^{X \times X}_{A \times A,B \times B}\left( c \circ \Tang( f ) \circ c \right): A \times A \to B \times B.
\]Asking that these maps be equal is equivalent to the other rules. 

\begin{proposition} A Cartesian differential category with a Conway operator is a traced Cartesian differential category if and only if for the induced trace operator $\Trace$ and for every map $f: A \times X \to B \times X$, the following equality holds:  
	\begin{align}\label{eq:trace-fix}
		\Tang(\Trace^X_{A,B}\left( f \right)  ) =\Trace^{\Tang X}_{\Tang A,\Tang B}( c \circ \Tang( f ) \circ c ) 
	\end{align}
which means that $(T,c)$ is a traced monoidal functor (\cite{joyal1996traced}) where the monoidal product is cartesian.
\end{proposition} 
\begin{proof} For the $\Rightarrow$ direction, let $f = \langle f_1, f_2 \rangle$. Let us first put $c \circ \Tang\left( f \right) \circ c$ into a pair so that we can take its trace. By \cite[Lemma 2.9(iv)]{lemay2021exponential}, recall that $\Tang(\langle h,k \rangle) = c \circ \langle\Tang(h),\Tang(k) \rangle$. From this it is straightforward to compute that 
	\[
	c \circ \Tang\left( f \right) \circ c = \left \langle\Tang(f_1) \circ c,\Tang(f_2) \circ c \right \rangle.
	\]Then, using the tangent-fixpoint rule, we compute that: 
	\begin{gather*} 
		\Trace^{\Tang X}_{\Tang A,\Tang B}( c \circ \Tang( f ) \circ c ) =\Trace^{\Tang X}_{\Tang A,\Tang B}(  \langle\Tang(f_1) \circ c,\Tang(f_2) \circ c  \rangle ) \\
		=\Tang(f_1) \circ c \circ  \langle 1_{A \times A},  \fix^{\Tang X}_{\Tang A}(\Tang(f_2) \circ c )  \rangle\\ 
		=\Tang(f_1) \circ c \circ  \langle\Tang(1_A),  \Tang( \fix^X_A(f_2) )   \rangle  \\
		=\Tang(f_1) \circ\Tang (  \langle 1_A,  \fix^X_A(f_1)  \rangle ) \\
		=\Tang( f_1  \circ  \langle 1_A,  \fix^X_A(f_2)  \rangle ) \\
		= \Tang(\Trace^X_{A,B}( \langle f_1, f_2 \rangle )  )  =\Tang(\Trace^X_{A,B}( f )  ) 
	\end{gather*}
	So (\ref{eq:trace-fix}) holds as desired. For the $\Leftarrow$ direction, we will show that (\ref{eq:trace-fix}) implies the tangent-fixpoint rule. So for $f: A \times X \to X$, we get: 
	\begin{gather*}
	\Tang ( \fix^X_A(f)) =  \Tang(\Trace^X_{A,X}( \langle f, f \rangle )  ) \\
	= \Trace^{\Tang X}_{\Tang A,\Tang X}\left( c \circ \Tang\left( \langle f, f \rangle \right) \circ c \right)\\
	 = \Trace^{\Tang X}_{\Tang  A,\Tang X}\left( \left\langle  \Tang\left( f \right), \Tang\left( f \right) \right \rangle \circ c \right) \\
	= \Trace^{\Tang X}_{\Tang A,\Tang X}\left( \left\langle  \Tang\left( f \right) \circ c, \Tang\left( f \right) \circ c \right \rangle \right) \\
	= \fix_{\Tang  A}^{\Tang X}\left(\Tang(f)\circ c\right)
	\end{gather*}
	So the tangent-fixpoint rule holds. Then by Proposition \ref{prop:tanfix-difffix}, it follows that we have a traced Cartesian differential category.  
\end{proof}

\section{Fixpoints and Linearity}\label{subsec:fix-linear}

An important notion in a Cartesian differential category is the concept of \emph{linearity}. In particular, linear maps are those whose derivatives are simply themselves evaluated in the second argument. The term linear is justified since in many models, linear in the Cartesian differential category sense corresponds to being linear in the classical algebraic sense. 

\begin{definition} In a Cartesian differential category, a map $f: A \to B$ is \textbf{linear} \cite[Definition 2.1.1]{blute2009cartesian} if $\mathsf{D}[f] = f \circ \pi_2$, which in the term calculus is expressed as: 
	\[\frac{\mathsf{d}f(x)}{\mathsf{d}x}(a) \cdot b = f(b)\] 
	For a Cartesian differential category $\cat{C}$, its subcategory of linear maps will be denoted as $\mathsf{LIN}[\cat{C}]$.
\end{definition}

Properties of linear maps can be found in \cite[Lemma 2.6]{cockett2020linearizing}, such as the fact that they are closed under composition, sum, and product structure. In particular, every linear map is also additive \cite[Lemma 2.6.(i)]{cockett2020linearizing}. From this, it follows that $\mathsf{LIN}[\cat{C}]$ has finite \emph{biproducts} \cite[Corollary 2.2.3]{blute2009cartesian}. It is important to note that although linear maps and additive maps often coincide, in an arbitrary Cartesian differential category not every additive map is necessarily linear. 

We now consider how linearity interacts with parametrized fixpoint operators. The first natural thing to ask is that the linear maps be the \emph{strict} maps. Indeed, an important concept for fixpoint operators is the notion of \emph{uniformity}, which is used to characterize fixpoint operators uniquely without relying on order-theoretic arguments. The uniformity axiom is then relative to a subcategory of \emph{strict} maps. 

\begin{definition}\label{def:strictmaps} In a category $\cat{C}$ with finite products and equipped with a parametrized fixpoint operator $\fix$, a map $h: X \to Y$ is \textbf{strict} \cite[Definition 4.4]{hasegawa2004uniformity} if for every $f:A \times X \to X$ and $g : A \times Y\to Y$, if the equality on the left holds, then the equality on the right holds:
	\[
	h \circ f = g \circ (1_A \times h) \quad \Rightarrow \quad h \circ \fix^X_A(f)= \fix^Y_A(g).
	\]
	which in the term calculus is expressed as:
	\[ h\left( f(a,x) \right) = g(a, h(x)) \quad \Rightarrow \quad h\left( \mu x. f(a,x) \right) = \mu y. g(a,y).  \]
\end{definition}

\begin{definition}\label{def:uniffix} Let $\cat{C}$ and $\cat{S}$ be categories with finite products, and $\J: \cat{S} \to \cat{C}$ be a bijective-on-objects and strict product preserving functor. Then a parametrized fixpoint operator on $\cat{C}$ is \textbf{$\J$-uniform} \cite[Definition 2.8]{simpson2000complete} if for every map $h$ in $\cat{S}$, $\J(h)$ is strict in $\cat{C}$. 
\end{definition}

Therefore, $\cat{S}$ is understood as a suitable subcategory of strict maps of $\cat{C}$. For categories of domain-like structures, the strict maps are usually the ones which preserve bottom elements. In categorical models of Linear Logic with fixpoints, the strict maps of interests are the maps of the base monoidal category. For Cartesian differential categories, it is natural to ask that the linear maps be strict. This can also be expressed in terms of the uniformity axiom by considering the canonical inclusion functor $\J: \mathsf{LIN}[\cat{C}] \to \cat{C}$. Therefore, we may ask that the parametrized fixpoint operator be $\J$-uniform, which we call being \emph{linearly uniform}. Explicitly: 

\begin{definition} A parametrized fixpoint operator on a Cartesian differential category is \textbf{linearly uniform} if every linear map is strict. 
\end{definition}

In Section \ref{sec:bifreefixdif}, we will use the uniformity axiom to make use of a theorem by Plotkin and Simpson \cite{simpson2000complete} so that we can prove the necessary compatibility results between differential combinators and parametrized fixpoint operator in coKleisli categories of comonads satisfying some conditions on bifree algebras. 

Linearity also behaves quite nicely with Conway operators. An important first observation is that for a Conway operator, the parametrized fixpoint of a linear map is again linear.

\begin{lemma} In a traced Cartesian differential category, if a map $f: A \times X \to X$ is linear, then $\fix^X_A(f): A \to X$ is also linear.
\end{lemma}
\begin{proof} We prove this using the term calculus. Using the strong differential-fixpoint rule and the fact that $f$ is linear, we compute: 
	\begin{gather*}
		\frac{\mathsf{d} \mu x. f(u,x)}{\mathsf{d}u}(a) \!\cdot\! b \!=\!\mu y. \frac{\mathsf{d}f(u,v)}{\mathsf{d}(u,v)} \left( a, \mu x. f(a,x) \right) \!\cdot\! \left( b, y \right)  = \mu y. f(b,y) 
	\end{gather*}
	and we conclude that $\fix^X_A(f)$ is also linear. 
\end{proof}

Therefore, if $\cat{C}$ is a traced Cartesian differential category with Conway operator $\fix$, then $\fix$ is also a Conway operator on $\mathsf{LIN}[\cat{C}]$. This in turn means that $\mathsf{LIN}[\cat{C}]$ also inherits the trace operator. However, for a category with finite biproducts, to give a trace operator is equivalent to giving a \emph{repetition operator} \cite[Proposition 6.11]{selinger2010survey}. So for a traced Cartesian differential category, its Conway operator induces a repetition operator on its subcategory of linear maps.

\begin{definition}\label{def:repetitionop}
	For a category $\cat{C}$ with finite biproducts, a \textbf{repetition operator} \cite[Proposition 6.11]{selinger2010survey} is a family of functions $(-)^*$ indexed by objects in $\cat{C}$, 
	\begin{align*}
		(-)^*_X: \cat{C}(X, X) \to \cat{C}(X,X) && \infer{f: X  \to X}{f^*: X \to X}
	\end{align*}
	such that the following equalities hold for all $f,g :X \to X$: \\
	\noindent 1. \textbf{Fixpoint axiom}: $f^* = 1_X + f\circ f^*$ \\ 
	\noindent 2. \textbf{Addition axiom}: $(f+g)^* = (f^*\circ g)^* \circ f^*$ \\ 
	\noindent 3. \textbf{Dinaturality axiom}: $(f\circ g)^* \circ f = f\circ (g\circ f)^*$
\end{definition}

\begin{corollary}\label{cor:linrep} For a traced Cartesian differential category $\mathbb{C}$, $\mathsf{LIN}[\cat{C}]$ has a repetition operator defined as $f^* = \fix^X_X\left( \pi_1 + f \circ \pi_2 \right)$, which in the term calculus is expressed as follows: 
	\begin{equation}\label{eqn:repetitionop}
		f^*(a) = \mu x.\left( a + f(x) \right).
	\end{equation}
\end{corollary}

On the other hand, we can also consider linearity in a certain argument. For the story of this paper, we are particularly interested in maps which are linear in their second argument.  

\begin{definition} In a Cartesian differential category, a map $f: A \times X \to B$ is \textbf{linear in its second argument} \cite[Definition 4.5]{cockett2020linearizing} (or simply linear in $X$ when there is no confusion) if \[
	\Diff[f] \circ \langle a,x,0,y \rangle = f \circ \langle a,y \rangle
	\]
	which in the term calculus is expressed as: 
	\[\frac{\mathsf{d}f(u,v)}{\mathsf{d}(u,v)}(a,x) \cdot (0,y) = f(a,y).\] 
\end{definition}

In particular, {\bf [CD.6]} is precisely the statement that for any map $f: A \to B$, its derivative $\Diff[f]: A \times A \to B$ is linear in its second argument $A$. Similarly, for a map ${f: A_1 \times \hdots \times A_n \to B}$, we can also define what it means to be linear in its $j$-th argument $A_j$ \cite[Definition 2.6]{garner2020cartesian}. Now, taking the parametrized fixpoint in a linear argument results in zero. This makes sense since a map $f: A \times X \to X$ which is linear in $X$ is also additive in $X$. In particular, $f(a,0) = 0$. Thus $0$ is a parametrized fixpoint. It is important to note that $f: A \times X \to X$ being linear in $X$ is quite different from $f$ being linear. Thus the following lemma does not clash with the previous one. 

\begin{lemma} In a traced Cartesian differential category, if a map $f: A \times X \to X$ is linear in $X$ then $\fix^X_A(f) = 0$.  
\end{lemma}
\begin{proof} 
	We use the term calculus again. The key axiom here is \textbf{[CD.2]}, which recall says that evaluating a derivative at zero in its second argument is zero. Therefore, using the strong differential-fixpoint rule, we compute that $\mu y. f(a, y)$ is equal to 
	\begin{gather*} 
	\mu y. \frac{\mathsf{d}f(u,v)}{\mathsf{d}(u,v)} \left( a, \mu x. f(a, x)\right) \cdot \left( 0, y \right) = \frac{\mathsf{d} \mu x. f(u,x)}{\mathsf{d}u}(a) \cdot 0 = 0 \qedhere
	\end{gather*}
\end{proof}

\section{Closed Setting}\label{sec:closed}

Let us now consider the case when the underlying category is Cartesian \emph{closed}. This is an important case since Cartesian \emph{closed} differential categories (also sometimes called differential $\lambda$-categories) provide the categorical semantics of the differential $\lambda$-calculus. For a more in-depth introduction to Cartesian closed differential categories and the differential $\lambda$-calculus, we refer the reader to~\cite{bucciarelli2010categorical,manzonetto2012categorical,Cockett-2019, cockett2020linearizing, ehrhard2003differential}. 

For a Cartesian \emph{closed} category, we denote the internal hom by $A \Rightarrow B$, the evaluation map by $\eval_{A,B}: \left( A \Rightarrow B \right) \times A \to B$, and for a map $f: C \times A \to B$, we denote its currying by $\lambda(f): C \to A \Rightarrow B$, which recall is the unique map such that $\eval_{A,B} \circ \left(\lambda(f) \times 1_A \right) = f$. As explained in \cite[Lemma 4.10]{Cockett-2019}, there are two equivalent ways of expressing compatibility between the closed structure and the differential combinator: one in terms of the Curry operator and one in terms of the evaluation map. In terms of the latter, a \textbf{Cartesian closed differential category} \cite[Section 4.6]{Cockett-2019} is a Cartesian differential category which is Cartesian closed such that every evaluation map is linear in their internal hom argument (which by our convention means they are linear in their \emph{first} argument). Another equivalent way of axiomatizing Cartesian closed differential categories is via an axiom which says that the derivative of a Curry is the Curry of the \emph{partial} derivative \cite[Definition 6.2]{cockett2020linearizing}. Every model of the differential $\lambda$-calculus induces a Cartesian closed differential category \cite[Theorem 4.3]{Cockett-2019}, and conversely, every Cartesian closed differential category induces a model of the differential $\lambda$-calculus \cite[Theorem 4.12]{bucciarelli2010categorical}.

On the other hand, in a Cartesian closed category, fixpoint operators can be internalized with the notion of \emph{fixpoint combinator}:

\begin{definition}\label{FixpointCombinator} For a Cartesian closed category $\cat{C}$, a \textbf{fixpoint combinator} is a family of maps $\{\Y_X :X \Rightarrow X \to X\}_X$ indexed by the objects of $\cat{C}$, satisfying the following axiom for all $X$:
	\[\Y_X = \eval \circ \langle 1_{X \Rightarrow X}, \Y_X \rangle.\]
\end{definition}
The equivalence between parametrized fixpoint operators and fixpoint combinators in the Cartesian closed setting works as follows: starting with a parametrized fixpoint operator $\fix$, for every object $X$, define the map $\Y_X: X \Rightarrow X \to X$ as the parametrized fixpoint of the evaluation map $\eval_{X,X}: \left( X \Rightarrow X \right) \times X \to X$, that is,
\[
\Y_{\scriptsize X} := \fix^X_{\scriptsize X \Rightarrow X}( \eval_{X,X}).
\] For the other direction, for every $f: A \times X \to X$, its parametrized fixpoint can be expressed in terms of $\Y_X$ as: 
\begin{align}\label{eq:fixmap}
	\fix^X_A(f) = \Y_X \circ \lambda(f).
\end{align}
From this identity, we can express the differential-fixpoint rule and the others in terms of $\Y$. In particular, we highlight that (\ref{ax:fixtan-closed}) is the analogue of the equation Ehrhard computes in the coKleisli category of a \emph{coherent} differential category which is Scott \cite[Theorem 5.29]{ehrhard2023coherent}. 

\begin{proposition} In a Cartesian closed differential category with a parametrized fixpoint $\fix$,
	\begin{enumerate}[{\em (i)}]
		\item $\fix$ satisfies the differential-fixpoint rule if and only if the following equality holds: 
		\begin{align}\label{ax:fixdiff-closed}
			\Diff \left[  \Y_X \right] = \pi_2 \circ \Y_{\Tang X} \circ \lambda\left( \Tang(\eval_{X,X}) \circ c \right) 
		\end{align}
		\item $\fix$ satisfies the tangent-fixpoint rule if and only if the following equality holds: 
		\begin{align}\label{ax:fixtan-closed}
			\Tang(  \Y_X ) =  \Y_{\Tang X} \circ \lambda\left( \Tang(\eval_{X,X}) \circ c \right)
		\end{align}
		\item $\fix$ satisfies the strong differential-fixpoint rule if and only if the following equality holds: 
		\begin{align}
			\Diff\!\left[ \Y_X \right] \!= \fix_X \circ \lambda\left( \Diff[\eval] \circ \left \langle \pi_1, \Y_X \circ \pi_1, \pi_2, \pi_3 \right \rangle \right) 
		\end{align}	
	\end{enumerate}
\end{proposition}
\begin{proof} Let us prove (i). For the $\Rightarrow$ direction, by setting $f = \eval_{X,X}$ in (\ref{ax:fixdiff}), and then using (\ref{eq:fixmap}), we get that: 
	\[\begin{aligned}
		\Diff \left[  \Y_X \right] &= \Diff [  \fix^X_{X \Rightarrow X}( \eval_{X,X} )] \\
	&=  \pi_2 \circ \fix_{\Tang X \Rightarrow \Tang X}^{\Tang X}
	\left( \Tang(\eval_{X,X}) \circ c \right) \\ 
	&= \pi_2 \circ \Y_{\Tang X} \circ \lambda\left( \Tang(\eval_{X,X}) \circ c \right)
	\end{aligned}\]
	For the $\Leftarrow$ direction, first recall that $\lambda(g) \circ h = \lambda\left( g \circ (h \times 1) \right)$. We also have that $\Tang$ is a strong Cartesian monoidal functor \cite[Lemma 2.9(v)]{lemay2021exponential}
	\[c \circ \left(\Tang(f) \times \Tang(g) \right)  =  \Tang(f \times g) \circ c.\]From these identities, it follows that \[
	\lambda\left( \Tang(f) \circ c \right) \circ \Tang(g) = \lambda\left( \Tang\left(f \circ (g\times 1) \right) \circ c \right).
	\] 
	Using this and (\ref{eq:fixmap}), we compute that: 
	\[\begin{aligned}
		\Diff[\fix^X_A(f)] &= \Diff \left[ \Y_X \circ \lambda(f)  \right] \\
		&= \Diff \left[ \Y_X   \right] \circ \Tang\left(\lambda(f) \right) \\
		&= \pi_2 \circ \Y_{\Tang X } \circ \lambda\left( \Tang(\eval_{X,X}) \circ c \right) \circ \Tang\left(\lambda(f) \right) \\
		&= \pi_2 \circ \Y_{\Tang X} \circ \lambda\left( \Tang\left(\eval_{X,X} \circ \left(\lambda(f)  \times 1 \right) \right)  \circ c \right)  \\
		&= \pi_2 \circ \Y_{\Tang X} \circ \lambda\left( \Tang\left( f \right)  \circ c \right) \\
		&= \pi_2 \circ \fix_{\Tang A}^{\Tang X}\left(\Tang(f)\circ c \right)
	\end{aligned}\]
	So the differential-fixpoint rule holds as desired. One can prove (ii) and (iii) using similar computations. 
\end{proof}

\section{Examples}\label{sec:examplesfixdif}

\subsection{Categories with finite biproducts}\label{sec:biproduct}
Any category $\cat{C}$ with finite biproducts is a Cartesian differential category where the differential combinator is given by pre-composition with the second projection: $\Diff[f] = f \circ \pi_2$ \cite[Example 2.10]{cockett2020linearizing}. If $\cat{C}$ comes equipped with a parametrized fixpoint operator such that the projection maps are strict (Definition \ref{def:strictmaps}), then the differential-fixpoint rule holds, and so $\cat{C}$ is a Cartesian differential fixpoint category. On the other hand, any Conway operator on $\cat{C}$ can easily be seen to satisfy the strong differential-fixpoint rule. Therefore, any category with finite biproducts which is equipped with a Conway operator (or equivalently a trace operator or a repetition operator) is a traced Cartesian differential category. 
\subsection{Categories of (weighted) relations}\label{ex:relations}
The category $\Rel$ whose objects are sets and morphisms are binary relations between them is one of the most basic models of linear logic, with many other models arising as refinements of
it. The operation mapping a set $A$ to the set of finite multisets over $A$, \[
\oc A := \{ {m : A \to \NN} \mid m \text{ has finite support}\}
\]
can be equipped with a comonad structure on $\Rel$. The induced co-Kleisli category $\Rel_\oc$ is a Cartesian differential category \cite[Section 5.1]{bucciarelli2010categorical} and also has a canonical Conway operator \cite[Proposition 3]{grellois2015infinitary}.

Binary relations over sets can be generalized to weighted relations over a continuous semi-ring \cite{laird2013weighted}.
A \emph{continuous semi-ring} is a semi-ring $(\Sring, \leq,  \Zero,+, \One, \cdot)$ which is also a cpo, whose zero element $\Zero$ is the bottom, and both the addition and multiplication are continuous. For a continuous semi-ring $\Sring$, the category $\Rel^{\Sring}$ has objects sets and maps from $A$ to $B$ are functions $R: A \times B \to \Sring$. The composite of $R : A \times B \to \Sring$ and $S: B \times C \to \Sring$ is the function $S \circ R: A \times C \to \Sring$ given by 
\[(S \circ R)(a,c) = \sum_{b\in B} R(a,b) \cdot S(b,c)\] where the sum on the right is well defined by the cpo structure of $\Sring$. The finite multiset construction also induces a comonad on $\Rel^{\Sring}$ and the corresponding co-Kleisli category is a Cartesian differential category. The derivative of a weighted relation $R: \oc A \times B \to \Sring $ is $\Diff[R] : \oc A \times \oc A \times B \to \Sring$:
\[
(m, n, b) \mapsto \begin{cases}
	R(m + n, b) & \text{ if } n = [a] \textit{ for some } a \in A \\
	0 & \text{otherwise}
\end{cases}
\]
The co-Kleisli category $\Rel^{\Sring}_\oc$ also has a Conway operator obtained via Kleene iteration. For a weighted relation $R: \oc A \times \oc X \times X \to \Sring$, its least fixpoint is $\fix R : \oc A \times X \to \Sring$ inductively defined as $\fix R := \bigvee_{n\in \NN} \fix_n R$, where for $m \in \oc A$, we define 
\[
\fix_0 R (m) = 0 \qquad \fix_{n+1} R (m) = R(m, \fix_n R (m))
\]
with composition in this expression being taken in the co-Kleisli category. The fixpoint operator defined above satisfies the universal property of being the unique parametrized fixpoint operator that is uniform with respect to linear morphisms, which here correspond precisely to the morphisms of $\Rel^{\Sring}$. This uniform Conway operator and the differential combinator verify the (strong) differential-fixpoint rule, and therefore $\Rel^\Sring_{\oc}$ is a traced Cartesian closed differential category, and moreover the Conway operator is linearly uniform. This is true because this example is an instance of a more general statement concerning fixpoint operators obtained via bifree algebras of comonads \cite{simpson2000complete} -- which we discuss in Section \ref{sec:bifreefixdif} below. When taking the Boolean semi-ring $\mathbb{B} = \{\true, \false\}$, we get back the relational model, $\Rel^{\mathbb{B}}_\oc = \Rel_\oc$.

\subsection{Formal Power Series}

For a continuous semi-ring $\Sring$, we may consider the full subcategory of $\Rel^{\Sring}_\oc$ whose objects are the finite sets. This subcategory of finite sets is again a traced Cartesian differential category and is equivalent to the Lawvere theory of formal power series over $\Sring$. We obtain a traced Cartesian differential category of formal power series over any continuous semi-ring. Even though this is a subcategory of the previous example, since formal power series over continuous semi-rings are used in many areas of computer science, it is worthwhile to spell out some details. 

For a continuous semi-ring $\Sring$, let $\Pow^{\Sring}$ be the category whose objects are the natural numbers $n \in \NN$ and where a map $p: n \to m$ is an  $m$-tuple of power series in $n$-variables:
\[
p = (p_1(x_1, \dots, x_n),\dots, p_m(x_1, \dots, x_n)).
\] Since coefficients are over a continuous semi-ring, the composition of power series is well-defined. Now for $n\in \NN$, define the set $\ints{n}$ as $\{1, \dots, n\}$. Then the full subcategory $\Rel^{\Sring}_\oc$ with objects restricted to finite sets of the form $\ints{n}$ is isomorphic to $\Pow^{\Sring}$. Indeed, a weighted relation $R : \oc \ints{n}  \times \ints{m} \to \Sring$ in $\Rel^{\Sring}_\oc$ corresponds to the tuple of power series $p = (p_1 (\vec{x}), \dots, p_m(\vec{x}))$ where for $1\leq j \leq m$: 
\[
p_j (x_1, \dots, x_n) = \sum\limits_{(k_1, \dots, k_n) \in \Sring^n} R(k_1, \dots, k_n, j) \cdot x_1^{k_1} \dots x_n^{k_n}
\]
Therefore, it follows that $\Pow^{\Sring}$ is a traced Cartesian differential category. In particular, the least fixpoint in the weighted relational model corresponds to the least fixpoint for power series over continuous semi-rings first described in \cite{rozenberg2012handbook}. On the other hand, for a map $p: \ints{n} \to 1$, which is just a power series in $n$ variables $p(x_1, \dots, x_n)$, applying the differential combinator to it results in the sum of its partial derivatives: 
\[
\Diff[p](\vec{x}, \vec{a}) = \sum_{i=1}^n \frac{\partial p(\vec{x})}{\partial x_i}\cdot a_i
\]
This corresponds to the notation $\Diff p \vert_{\vec{x}}(\vec{a})$ in papers such as \cite{esparza2010newtonian}. By \textbf{[CD.4]}, for a tuple $p = (p_1 (\vec{x}), \dots, p_m(\vec{x}))$, its derivative is \[\Diff[p] = \left(\Diff[p_1](\vec{x}, \vec{a}), \dots, \Diff[p_m](\vec{x}, \vec{a}) \right).\]
\subsection{Quantale-enriched profunctors}\label{ex:quantale}
Another possible generalization of the notion of relations is given by the notion of quantale enriched profunctors \cite{kelly1982basic, Stubbe2004}. A \emph{quantale} $\QQ$ is a complete lattice that is symmetric monoidal closed such that the tensor distributes over arbitrary suprema. Quantales are idempotent semi-rings with join $\vee$ as additive structure and tensor $\otimes$ as multiplicative structure. Instead of having sets as objects, we can consider richer structures given by the notion of quantale enriched categories: for a quantale $\QQ$, a \emph{$\QQ$-category} $\cat{A}$ consists of a set of objects $\Ob (\cat{A})$ and for all objects $a,b$ an object $\cat{A}(a,b)$ in $\QQ$ together with composition $\cat{A}(a,b) \otimes \cat{A}(b,c) \leq \cat{A}(a,c)$ and identity inequalities $\One \leq \cat{A}(a,a)$ in $\QQ$. For example, a $\QQ$-category for the Lawvere quantale $\Rp\cup \{ \infty\}$ corresponds to the notion of \emph{generalized metric space}~\cite{lawvere1973metric}.

For $\QQ$-categories $\cat{A}, \cat{B}$, a $\QQ$-profunctor $R: \cat{A} \proft \cat{B}$ (also called distributors or bimodules) from $\cat{A}$ to $\cat{B}$ is a function $R : \Ob (\cat{A}) \times \Ob(\cat{B}) \to \QQ$ with biaction inequalities for all $a, a' \in \Ob(\cat{A})$, $b,b' \in \Ob(\cat{B})$: 
\[
\cat{A}(a,a') \otimes R(a,b) \otimes \cat{B}(b',b) \leq R(a',b').
\]There is an analogous free exponential construction mapping a $\QQ$-category $\cat{A}$ to the $\QQ$-category $\oc \cat{A}$ whose set of objects is the set of finite sequences $\vec{a} = \lis{a_1, \dots, a_n}$ of objects of $\cat{A}$ and for $\vec{a} = \lis{a_1, \dots, a_n}$ and $\vec{b} = \lis{b_1, \dots, b_m}$, the hom-object is given by
\[
\oc \cat{A}(\vec{a}, \vec{b}) = \begin{cases}
	\bigvee\limits_{\sigma : n  \cong n} \bigotimes\limits_{1\leq i \leq n} \cat{A}(a_i, b_{\sigma{i}}) & \text{ if } n=m\\
	\bot & \textit{ otherwise}
\end{cases}
\]
The induced co-Kleisli category has a differential combinator, where the derivative of a $\QQ$-profunctor $R : \oc \cat{A} \proft \cat{B}$ generalizes the one in the discrete case for weighted relations:
\[
(\lis{a_1, \dots, a_n},\! \lis{a_1', \dots, a_m'}, b) \!\mapsto\!\! \begin{cases}
	R(\lis{a_1, \dots, a_n, a_1'}, b) & \text{ if } m = 1\\
	\bot & \text{otherwise}
\end{cases}
\]
The induced co-Kleisli category also has a fixpoint operator mapping a $\QQ$-profunctor $ R : \oc \cat{A} \otimes \oc \cat{X} \proft \cat{X}$ to $\fix R : \oc \cat{A} \proft \cat{X}$ obtained by Kleene iteration as $\vee_{n \in \NN} \fix_n R$ with $\fix_0 R =\bot$ and 
\[
\fix_{n+1}R (\vec{a}, x) = \bigvee\limits_{\substack{\vec{y}= \lis{y_1, \dots, y_n} \in \oc \cat{X} \\\vec{a}_0, \dots, \vec{a}_n \in \oc \cat{A}}} R(\vec{a}_0, \vec{y}, x) \otimes\bigotimes\limits_{1\leq i \leq n} \fix_n R(\vec{a}_i, y_i)  
\]
So we get a Cartesian differential fixpoint category, which is again an instance of the bifree algebras story in the next section.

\subsection{Fixpoints from bifree algebras}\label{sec:bifreefixdif}

Recall that for an endofunctor $F : \cat{C} \to \cat{C}$, a \emph{bifree $F$-algebra} \cite[Section 5]{simpson2000complete} is an initial $F$-algebra $(A, f: FA \to A)$ such that the inverse of $f$ is a final $F$-coalgebra $(A, f^{-1}: A \to FA)$. A result by Simpson and Plotkin allows one to construct a parametrized fixpoint operator in the co-Kleisli category of a comonad whose underlying endofunctor has suitable bifree algebras~\cite{simpson2000complete}. We proceed to show that if the co-Kleisli category is also a Cartesian differential category, then it is also a Cartesian differential fixpoint category. 

\begin{theorem}\cite[Proposition 6.5 \& Theorem 3]{simpson2000complete}\label{thm:PlotkinSimpson}
	Let $\cat{C}$ be a category with finite products equipped with a comonad $(F, \delta, \varepsilon)$. Let $\J : \cat{C} \to \cat{C}_{F}$ be the free functor induced by the comonadic adjunction. 
	\begin{enumerate}
		\item If for all objects $A$ in $\cat{C}$, the endofunctors $F( A \times -)$ have a bifree algebra, then $\cat{C}_{F}$ has a unique uniform (with respect to $\J$) parametrized fixpoint operator.
		\item If for all objects $A$ in $\cat{C}$, the endofunctors $F( A \times F(A \times -))$ and $F(A \times - \times -)$ have a bifree algebra, then $\cat{C}_{F}$ has a unique uniform (with respect to $\J$) Conway fixpoint operator.
	\end{enumerate}
\end{theorem}

We give a sketch of how the parametrized fixpoint operator is constructed via bifree algebras. Assume that for all $A$ in $\cat{C}$, the endofunctor $F( A \times -)$ has a bifree algebra $w_A : F(A \times \Phi_A) \to \Phi_A$ which we can view as a map $A \times \Phi_A \to \Phi_A$ in the co-Kleisli category $\cat{C}_{F}$. As explained in \cite[Section 6]{simpson2000complete}, for all $A$ in $\cat{C}$, there exists a unique map $t_A : A \to \Phi_A$ in $\cat{C}_{F}$ such that
\[
w_A \circ (1_A \times t_A) \circ \Delta_A =  t_A.
\] 
Moreover, for all $f : A \times X \to X$ in $\cat{C}_{F}$, there exists a unique map $u_f : \Phi_A \to X$ in $\cat{C}$ such that 
\[
\J(u_f) \circ w_A = f \circ \J(1_A \times u_f).
\] 
Then, for a map $f : A \times X \to X$ in $\cat{C}_{F}$, its parametrized fixpoint is defined as the following composite: 
\[
\fix^X_A(f) := A \xrightarrow{t_A} \Phi_A \xrightarrow{\J(u_f)} X.
\]
Now suppose that $\cat{C}_{F}$ is also a Cartesian differential category, such that for all maps $f$ in $\cat{C}$, $\J(f)$ is linear in $\cat{C}_{F}$. We will now show that the parametrized fixpoint operator defined above satisfies the tangent-fixpoint rule. 

To do so, the key to this proof is that if $k$ is linear, then $\Tang(k) = k \times k$ \cite[Lemma 2.(iii)]{lemay2021exponential}. In particular, for all $g$ in $\cat{C}$, we have that $\Tang(\J(g)) = \J(g) \times \J(g)$. Now, to help with notation, define $e : \Phi_{A \times A} \to \Phi_A \times \Phi_A$ as $e := u_{\Tang(w_A) \circ c}$, that is, the unique map such that 
\[
\J(e) \circ w_{A\times A} = \Tang(w_A) \circ c \circ (1_{A \times A}\times \J(e)).
\]
Then, for $f : A \times X \to X \in \cat{C}_{F}$, we compute that (note $\J$ also preserves products strictly): 
\[\begin{aligned}
	&\Tang(f) \circ c \circ ( 1_{A \times A} \times \J((u_f \times u_f) \circ e ) ) \\ 
	&=\Tang(f) \circ c \circ (1_A \times 1_A \times \J(u_f)\times \J(u_f)) \circ \left(1_{A \times A}\times \J(e)\right) \\
	&= \Tang(f) \circ ( 1_A \times \J(u_f) \times 1_A \times  \J(u_f)) \circ c \circ \left(1_{A \times A} \times \J(e)\right) \\
	&= \Tang(f) \circ ( \J(1_A \times u_f) \times \J(1_A \times u_f))  \circ c \circ \left(1_{A \times A} \times \J(e)\right) \\ 
	&= \Tang(f) \circ \Tang ( \J(1_A \times u_f))  \circ c \circ \left(1_{A \times A} \times \J(e)\right) \\ 
	&= \Tang(\J(u_f)) \circ \Tang ( w_A )  \circ c \circ \left(1_{A \times A} \times \J(e)\right) \\ 
	&= (\J(u_f) \times \J(u_f)) \circ \J(e) \circ w_{A\times A}\\
	& = \J((u_f \times u_f) \circ e) \circ w_{A\times A}
\end{aligned}\]
By uniqueness, it follows that $u_{\Tang(f)\circ c} = (u_f \times u_f) \circ e$. Dually, since $w_{A\times A}^{-1}$ is a final coalgebra, there is a unique map $p: \Phi_A \times \Phi_A \to \Phi_{A \times A}$ in $\cat{C}$ such that 
\[
\J(p) \circ \Tang(w_A) \circ c = w_{A\times A} \circ (1_{A \times A} \times \J(p))
\] and also such that $\J(p) \circ \Tang(t_A) = t_{A\times A}$, which by uniqueness implies that $e\circ p = 1_{\Phi_A \times \Phi_A}$. Using these identities, we can finally compute that: 
\[\begin{aligned}
	&\fix_{\Tang A}^{\Tang X}(\Tang(f)\circ c)=\J(u_{\Tang(f) \circ c}) \circ t_{A \times A} \\
	&= \J((u_f \times u_f) \circ e) \circ \J(p) \circ \Tang(t_A) \\
	&= (\J(u_f) \times \J(u_f)) \circ \J(e) \circ \J(p) \circ  \Tang(t_A)\\ 
	&= (\J(u_f) \times \J(u_f))  \circ  \Tang(t_A) \\
	&= \Tang( \J(u_f) \circ t_A)= \Tang ( \fix^X_A(f))
\end{aligned}\]
We conclude that the tangent-fixpoint rule holds and that therefore $\cat{C}_{F}$ is a Cartesian differential fixpoint category. This construction subsumes both the example of weighted relations (Section \ref{ex:relations}) and the example of quantale profunctors (Section \ref{ex:quantale}) since the free exponential comonad in both examples has the required bifree algebras obtained via enrichment arguments. Both categories of weighted relations and $\QQ$-profunctors are cpo-algebraically compact~\cite{smyth1982category, cpo1992remarks,fiore2004axiomatic, freyd2006algebraically}: it means that every endofunctor that is cpo-enriched has a bifree algebra and one can verify that the endofunctor of the free exponential comonad in both cases is a cpo-enriched functor. For the weighted relational model, Laird computes these bifree algebras explicitly using nested finite multisets in \cite{laird2016fixed}.

\subsection{Fixpoint operators from fixpoint objects}\label{sec:lawvere}
A standard construction to obtain fixpoint operators in a Cartesian closed category is via cpo-enrichment. A category $\cat{C}$ is $\omega$-cpo enriched if each hom-set $\cat{C}$ is an $\omega$-cpo and composition $\cat{C}(A,B) \times \cat{C}(B, C) \to \cat{C}(A,C)$ is Scott-continuous. If in addition, the evaluation and pairing are monotonous, each hom-set $\cat{C}(A,B)$ has a bottom element $0$ and for all $f$, we have $0 \circ f = 0$ and $\eval \circ \langle 0, f \rangle = 0$, it was shown by Berry~\cite{berry1979modeles} that the category has a least parametrized fixpoint operator. For each $X$, we can construct a fixpoint combinator morphism $\Y_X: X \Rightarrow X \to X$ as the supremum $\Y_X = \bigvee_n \Y^n_X$ where \[
\Y_X^0 = 0 \quad \text{and} \quad \Y^{n+1}_X= \eval_{X,X} \circ \lis{1_{X\Rightarrow X}, \Y^n_X}.
\]
Explicitly, the parametrized fixpoint operator is defined as in (\ref{eq:fixmap}). Now, if $\cat{C}$ is a Cartesian differential category and the cpo bottom elements $0$ coincide with the zero of the additive structure, then we can show that this parametrized fixpoint operator satisfies the tangent-fixpoint rule, or equivalently $\Y_X$ satisfies (\ref{ax:fixtan-closed}). We omit the proof as it follows the same reasoning as done by Ehrhard in \cite[Theorem 5.29]{ehrhard2003differential} in coherent differential categories with a Scott fixpoint \cite{ehrhard2023coherent}. The known models of (coherent) PCF typically have reflexive objects (retractions $(X \Rightarrow X) \triangleleft X)$ which provide models of untyped $\lambda$-calculs and are particular cases of fixpoint objects. We prove below that the tangent fixpoint axiom holds for Cartesian closed categories where fixpoint objects induce a fixpoint operator by Lawvere's theorem~\cite{lawvere1969diagonal}.

\begin{theorem}\cite[Section 1]{lawvere1969diagonal}
	In a Cartesian closed category $\cat{C}$, if there is a morphism $r: X \to X \Rightarrow X$ such that for every morphism $q : A \to X \Rightarrow X$, there exists a morphism $ u : A \to X$ such that $q$ can be factored as $q = r \circ u$, then every morphism $f : A \times X \to X$ has a parametrized fixpoint $\fix^X_A(f): A \to X$.
\end{theorem}

The construction works as follows: for $f : A \times X \to X$, we first define the map $p :A \times X \to X$ as
\[
p:= f \circ (1_A \times \eval) \circ (1_A \times r \times 1_X) \circ (1_A \times \Delta_X).
\]
Now let $u : A \to X$ be a map factoring $\lambda(p):A\to X\Rightarrow X$ as $\lambda(p)= r \circ u$. Then the parametrized fixpoint is defined as follows: 
\[
\fix^X_A(f)=A \xrightarrow{u} X \xrightarrow{\Delta_X} X \times X \xrightarrow{r \times 1_X} X \Rightarrow X \times X \xrightarrow{\eval_{X,X}}X
\]

To prove the differential-fixpoint rule, we need to impose stronger conditions. So we assume that for every object $X$, there is a map $r_X : X \to X \Rightarrow X$ such that for every map $ q : A \to X \Rightarrow X$, there is a \emph{unique} morphism $u : A \to X$ such that $q = r_X \circ u$. This stronger assumption holds automatically for the retraction case with reflexive objects. We also need the choice of fixpoint objects to be compatible with the differential structure, so we assume that for all $X$,
\[
r_{\Tang X} = (\pi_1 \Rightarrow 1_{\Tang X }) \circ \Tang(r_X).
\]
Moreover, since $\pi_1 \circ \langle 1_X, 0 \rangle = 1_X$, we have that \[
\Tang(r_X) = (\lis{1_X, 0}  \Rightarrow 1_{\Tang X}) \circ r_{\Tang X}.
\] 
Now for a map $f : A \times X \to X$, define $s : \Tang A \times \Tang X \to \Tang X$ as the following composite: 
\begin{align*}
	s:=\Tang f \circ c \circ (1_{\Tang A} \times \eval)\circ (1_{\Tang A} \times r_{\Tang X} \times 1_{\Tang X}) \circ (1_{\Tang A} \times \Delta_{\Tang X}).
\end{align*}
Now let $w : \Tang A \to \Tang X$ be the unique map such that $r_{\Tang X} \circ w = \lambda (s)$. Note that $\Tang(\Delta_X) = c \circ \Delta_{\Tang X}$ so we obtain:  
\[
	\fix^{\Tang X}_{\Tang  A}(\Tang(f) \circ c) =  \eval_{\Tang X, \Tang X} \circ ( r_{\Tang X} \times 1_{\Tang X} ) \circ \Delta_{\Tang X} \circ w 
\]
On the other hand, to work out $\Tang(\fix^X_A(f))$, we use some identities that hold in a Cartesian closed tangent category \cite[Section 5.2]{gallagher2018differential}, so in particular hold in a Cartesian closed differential category. We first have that 
\[\Tang(\eval_{X,Y})  \circ c = \eval_{\Tang X, \Tang Y} \circ ( (\pi_1 \Rightarrow 1_{\Tang Y}) \times 1_{\Tang X } )\]
which allows us to first compute that $\Tang(p) \circ  c $ is equal to
\begin{gather*}
	\Tang(f \circ (1_A \times \eval) \circ (1_A \times r_X \times 1_X) \circ (1_A \times \Delta_X)) \circ c  \\
	= \Tang(f) \circ c \circ (1_{\Tang A} \times \eval) \circ (1_{\Tang A} \times r_{\Tang X} \times 1_{\Tang X}) \circ (1_{\Tang A} \times \Delta_{\Tang X})  =s
\end{gather*} 
Therefore, $w$ is the unique morphism such that $r_{\Tang X} \circ w = \lambda(\Tang(p) \circ c)$. Moreover, for a map $g : A \times B \to C$, we have that \[
\Tang(\lambda(g)) = (\lis{1_B, 0} \Rightarrow 1_{\Tang C}) \circ \lambda\left(\Tang(g) \circ c \right)
\] We then compute that: 
\begin{gather*}
	r_{\Tang X} \circ \Tang (u) = (\pi_1 \Rightarrow 1_{\Tang X}) \circ \Tang (r_X) \circ \Tang (u)\\
	=  (\pi_1 \Rightarrow 1_{\Tang X}) \circ \Tang (\lambda (p)) \\
	=  (\pi_1 \Rightarrow 1_{\Tang X}) \circ (\lis{1_X, 0}  \Rightarrow 1_{\Tang X}) \circ \lambda(\Tang (p) \circ c)\\
	=  \lambda(\Tang (p) \circ c)
\end{gather*}
By uniqueness, it then follows that $\Tang(u) = w$, so we may finally compute that:
\begin{gather*}
  \Tang(\fix^X_A(f)) = \Tang(\eval) \circ \Tang(r_X \times 1_X) \circ \Tang(\Delta_X) \circ \Tang(u) \\ 
	=  \Tang(\eval) \circ \Tang(r_X \times 1_X) \circ c \circ \Delta_{\Tang(X)} \circ w \\
	=  \Tang(\eval) \circ c \circ ( \Tang(r_X) \times 1_{\Tang X} ) \circ \Delta_{\Tang X} \circ w \\
	= \eval\circ ((\pi_1 \Rightarrow 1_{\Tang X}) \times 1_{\Tang X}) \circ ( \Tang(r_X) \times 1_{\Tang X} ) \circ \Delta_{\Tang X} \circ w \\
	= \eval \circ ( r_{\Tang X} \times 1_{\Tang X} ) \circ \Delta_{\Tang X} \circ w \\
= \fix^X_A(\Tang(f) \circ c)
\end{gather*}
So the tangent-fixpoint rule holds as desired. Note that the uniqueness requirement can be weakened if we assume that the factoring morphisms for $r_X$ and $r_{\Tang(X)}$ are chosen in a uniform way.

\section{Newton-Raphson iteration scheme}\label{sec:Newton}

	Newton-Raphson iteration has been extended to Kleene algebras or more generally power series over semi-rings~\cite{hopkins1999parikh, kiefer2007convergence, esparza2010newtonian, esparza2015fpsolve} where systems of equations in this setting can represent context-free grammars, data-flow equations, authorization problems, datalog queries etc. Newton iteration always converges for power series over $\omega$-continuous semi-rings and if we restrict to idempotent semi-rings, then it was shown to converge after a finite number of steps~\cite{esparza2010newtonian}. Similar ideas were also developed for enumerative combinatorics to compute efficiently large combinatorial structures that are defined via fixpoint equations~\cite{decoste1982approche, pivoteau2008generation, pivoteau2012algorithms}. In the combinatorial setting, the series always have positive coefficients corresponding to the number of structures of a given size and it also leads to a Newton iteration scheme that is always convergent.

	Since combinatorial species are first order terms of a Cartesian differential fixpoint (bi)category \cite{fiore2008cartesian} and power series over a continuous semi-ring are a subcategory of the weighted relational model, we were motivated to develop a general Newton-Raphson approximation scheme for Cartesian differential fixpoint categories which we present below. While the following construction works for the more general setting of Cartesian differential categories, we chose to present it in the Cartesian closed framework to emphasize the viewpoint of optimization procedures as higher order operators.
	\subsection{Newton approximants}
	We start by recalling the standard Newton-Raphson iteration method to find or approximate fixpoints for real-valued functions. For a differentiable function $f : \RR \to \RR$ verifying $f'(a) \neq 1$ for all $a \in \RR$, we consider a sequence of Newton approximants $\{z_n\}_{n \in \omega}$ with $z_0 \in \RR$ an initial chosen value and for each $n$, $z_{n+1} := \N(f)(z_n)$ where $\N(f) : \RR \to \RR$ is defined as
	\[
	\N(f) :  a \mapsto a + \frac{1}{1 - f'(a)} \cdot (f(a)-a).
	\]
	The sequence $\{z_n\}_{n \in \omega}$ may not converge or may not converge to a fixpoint of $f$ but under suitable conditions on $f$, it converges to a fixpoint and the convergence rate is quadratic. 
	Note first that the quotient $\frac{1}{1 - f'(a)} = 1 + f'(a) + (f'(a))^2 + \dots$ is a solution of the fixpoint equation 
	\[
	\frac{1}{1 - f'(a)} = 1 + f'(a) \cdot \frac{1}{1 - f'(a)}
	\]
	which corresponds exactly to the repetition operator fixpoint equation~\ref{def:repetitionop}. The general idea of Newton-Raphson iteration in our setting is to accelerate the computation of the non-linear fixpoint operator in the Cartesian (closed) category combining derivatives with the induced linear repetition operator in the category $\mathsf{LIN}[\cat{C}]$~(\ref{eqn:repetitionop}). 
	
	While we remain in a setting where everything is non-negative, we want to perform truncated subtractions: for example, in the semiring of non-negative extended reals $\RR_{\geq 0}^\infty=([0, \infty], \geq, \times, 1, +, 0)$, the truncated subtraction operation $a\ominus b$ corresponds to $\min\{0, a-b\}$. For a general Cartesian closed differential fixpoint category $\cat{C}$, we know that it is enriched over commutative monoids and we can equip each hom-set $\cat{C}(A,B)$ with the natural order relation induced by addition:
	 \[
	 \forall f,g : A \to B,  \; g \leq f \quad:\Leftrightarrow\quad  \exists h: A \to B,\; g+h=f.
	 \]
	  We assume that for every $f$ in $\cat{C}(A,B)$, the maps $f+(-)$ and $(-)+f$ have right adjoints for the $\geq$ ordering. Explicitly, for each $f,g,h$ in $\cat{C}(a,b)$, there exist maps $f \ominus g$ and $f \ominus h$ such that
	\[
	h \geq f \ominus g  \quad \Leftrightarrow \quad g + h  \geq f \quad \Leftrightarrow \quad  g \geq f \ominus h.
	\]
	It implies the following identities for all $f,g,h$ in $\cat{C}(a,b)$ and for every additive map $k : b\to c$:
	\begin{equation}
			\label{eqn:minusAdj}
		f \ominus (g + h) = (f \ominus g) \ominus h
	\end{equation}
	\begin{equation}
			\label{eqn:minusUnit}
		f \geq (f+g) \ominus g
	\end{equation}
	\begin{equation}
		\label{eqn:minusCounit}
		f \leq  (f\ominus g) +g \text{ and if }f \geq g, \text{ then }f = (f\ominus g) +g
	\end{equation}
	\begin{equation}
	\label{eqn:minusTriang}
	f \ominus g \leq (f \ominus h) + (h \ominus g)
\end{equation}
	\begin{equation}
		\label{eqn:minusAdd}
		k \circ (f \ominus g) \geq( k \circ f )\ominus (k \circ g)
	\end{equation}
	Since for $f : A \to A$ in $\cat{C}$, $\Diff[f] : A \times A \to A$ is linear in its second argument, it implies that for any point $\top \xrightarrow{a} A$,
	\[
	\Diff_a[f] :=\Diff[f]\circ (a\times 1_A) :A \to A
	\]
	is a morphism in $\mathsf{LIN}[\cat{C}]$. Therefore, $(\Diff_a[f])^* : A \to A$ verifies the fixpoint equation 
	\[
		(\Diff_a[f])^*\cdot b = b + 	\Diff_a[f] \cdot (	(\Diff_a[f])^* \cdot b).\] We can now define a map $\N : A \Rightarrow A \to A \Rightarrow A$ corresponding to Newton-Raphson iteration as follows:
	\[
	f \mapsto \lambda a. (a + (\Diff_a[f])^*\cdot (f(a) \ominus a))
	\]
	To define an analogue of Newton iteration in the parametrized Cartesian case, we recall first the notion of partial derivative. 
	For a map $f : A \times B \to C$, we can obtain its \emph{partial derivative}~\cite[Definition 6.2]{cockett2020linearizing} from the total derivative $\Diff[f] : A \times B \times A \times B \to C$ by precomposing with $0$ morphisms:
	\[\begin{aligned}
		\Diff_1[f] &:= \Diff[f] \circ (1_A \times 1_B \times \lis{1_A, 0}) : A \times B \times A \to C\\
		\Diff_2[f] &:= \Diff[f] \circ (1_A \times 1_B \times \lis{0, 1_B}) : A \times B \times B \to C
	\end{aligned}\]
	In term calculus, we write the partial derivatives as: 
	\begin{align*}
		\frac{\mathsf{d}f(x,b)}{\mathsf{d}x} (a) \cdot c = \frac{\mathsf{d}f(x,y)}{\mathsf{d}(x,y)} (a,b) \cdot (c,0) \\
		\frac{\mathsf{d}f(a,y)}{\mathsf{d}y} (b) \cdot c = \frac{\mathsf{d}f(x,y)}{\mathsf{d}(x,y)} (a,b) \cdot (0,c) 
	\end{align*}
	We can also define higher order derivatives by taking the partial derivative of the total derivative. For $f : A \to B$, we define the $n$-th derivative as $f^{(n)}  : =(\Diff_1)^n  [f] :A \times A^n \to B$, it is non-linear in its first argument and multi-linear in its last $n$ arguments. In the term calculus, we write this as:
	\[ f^{(n)}(a, (b, \dots, b)) = \frac{\mathsf{d}^nf(x)}{\mathsf{d} x^n}(a)\cdot (b, \dots, b) \]
	
	In the parametrized Cartesian case, we define an operator
	\[
	\N : \cat{C}(A \times X, X) \to \cat{C}(A \times X, X)
	\]
	mapping a morphism $f : A \times X \to X$ to the morphism
	\[
	(a,x) \mapsto x +\left(\frac{\mathsf{d}f(a,v)}{\mathsf{d}v}(x)\right)^*\cdot (f(a,x)\ominus x)
	\]
	To perform Newton iteration, we assume that $\cat{C}$ is a \emph{Taylor category}~\cite{ehrhard2018introduction, kerjean2023taylor,ehrhard2023coherent,ehrhard2023coherentbimonad}, which means that each map $f: A \to B$ is equal to its Taylor expansion:
\begin{equation}
	\label{eqn:Taylorsum}
	f(a+b) = \sum_n \frac{1}{n!} \frac{\mathsf{d}^nf(x)}{\mathsf{d}x^n}(a)\cdot (b, \dots, b)
\end{equation}
or equivalently for $b \leq a$
\begin{equation}
	\label{eqn:Taylorsubstract}
	f(a) =  \sum_n \frac{1}{n!} \frac{\mathsf{d}^nf(x)}{\mathsf{d}x^n}(b)\cdot (a\ominus b, \dots, a\ominus b)
\end{equation}

Assume that the fixpoint combinator $\Y : A \Rightarrow A \to A$ in $\cat{C}$ is obtained by Scott iteration as the supremum $\Y = \bigvee_n \Y_n$ with 
\[
\Y_0 =0 \quad \text{and}\quad \Y_{n+1} = \eval_{A,A} \circ \lis{1_{A \Rightarrow A}, \Y_n}
\]
where  $\cat{C}$ is assumed to be $\omega$-cpo enriched. We assume here that the $\omega$-cpo enrichment is compatible with the additive structure: $0$ is the bottom element and addition commutes with suprema of $\omega$-chains.

Since the fixpoint operator we consider is the least fixpoint operator, it implies that the induced linear repetition operator (Corollary \ref{cor:linrep}) satisfies the following induction axiom~\cite{esik2004inductive} for all $h, j, k \in  \mathsf{LIN}[\cat{C}](A,A)$: 
\begin{equation}
	\label{eqn:repetitionind}
	h + j \circ k \leq k \quad \Rightarrow \quad j^* \circ h \leq k.
\end{equation}
We can deduce it from the explicit computation of $j^* = \fix_A^A(\pi_1 + j \circ \pi_2)$ (Corollary \ref{cor:linrep}) which is equal to $\bigvee_{n} j^{*_{k}}$ with 
\[
j^{*_{0}}:=0\quad \text{and}\quad j^{*_{k+1}}:= 1_X + j \circ j^{*_k}.
\]
 We define a family of morphism $\Z_n : A \Rightarrow A \to A$ corresponding to the Newton approximants:
\[
\Z_0 =0 \quad \text{and} \quad \Z_{n+1} = \eval\circ \lis{\N, \Z_n}.
\]

	\begin{lemma}\label{lem:Newtonconv}
		For all $n$, the following inequalities hold:
		\[
		\Y_n \leq \Z_n \leq \Y \quad \text{and} \quad  \Z_n\leq \eval \circ \lis{1, \Z_n} \leq  \Z_{n+1} 
		\]
	\end{lemma}
	\begin{proof}
		For a fixed $f : A \to A$, to simplify the notation, we write $z_n:= \Z_n(f)$, $y_n := \Y_n(f)$, and $y := \Y(f)$,
		\begin{itemize}
			\item We have $z_0 = 0 \leq f(0)$ and we assume that for some $n$, $z_n \leq  f(z_n)$. Since $z_{n+1}$ is equal to
			\[
			\N(f)(z_n) = z_n + (\Diff_{z_n} [f])^* \cdot (f(z_n) \ominus z_n)
			\]
		we apply the Taylor formula (\ref{eqn:Taylorsum}) by setting $a := z_n$ and $b:=  (\Diff_{z_n} [f])^* \cdot (f(z_n) \ominus z_n)$, and obtain
			\[
			f(z_{n+1}) = f(a +b) \geq f(a) + \Diff_a[f] \cdot b
			= z_{n+1} \]
			where the last inequality is obtained from unfolding the repetition operator in the expression of $b$ and using the induction hypothesis to obtain that $f(z_n)= (f(z_n)\ominus z_n)+z_n$ by (\ref{eqn:minusCounit}).
			\item Using again the fixpoint axiom for the repetition operator and the fact that $z_n \leq f(z_n)$, we obtain that $z_{n+1}$ is equal to
			\[
				f(z_n) + (\Diff_{z_n}[f])\cdot(\Diff_{z_n}[f])^* \cdot (f(z_n)) \ominus x_n) \geq f(z_n).
			\]
			\item We show that for all $n$, $y_n \leq z_n$. The base case is trivial, and if we assume that $y_n \leq z_n$ for some $n$, then we have:
			\[
			y_{n+1} = f(y_n)\leq f(z_n)\leq z_{n+1}
			\]
			where the last two inequalities were proved above.
			\item The base case $z_0 =0 \leq y$ is trivial, we want to show that if for some $n$, $z_n \leq y$, then $z_{n+1}\leq y$. We show that for all $k \in \N$, 
			\[
			z_n + (\Diff_{z_n}[f])^{*_k} \cdot (f(z_n)) \ominus z_n) \leq f^{k+1}(z_n) \leq y
			\]
			which will imply the desired result by taking the supremum over $k$. For $k =0$, the inequality follows from $f(z_n) \geq z_n$. For the inductive step, we use the Taylor formula to obtain
			\begin{equation*}
			\begin{split}
				f^{k+2}(z_n) & \geq f (	z_n + (\Diff_{z_n}[f])^{*_k} \cdot (f(z_n)) \ominus z_n)) \\
				& \geq f(z_n) + (\Diff_{z_n}[f])\cdot(\Diff_{z_n}[f])^{*_k} \cdot (f(z_n)) \ominus z_n)) \\
				&= 	z_n + (\Diff_{z_n}[f])^{*_{k+1}} \cdot (f(z_n)) \ominus z_n)\qedhere 
			\end{split}\qedhere
		\end{equation*}
		\end{itemize}

	\end{proof}
	As a corollary of Lemma \ref{lem:Newtonconv}, we obtain that the Newton chain converges and approximates the least fixpoint from below as desired.

	\subsection{Taylor metric and convergence rate}
	To measure the convergence rate of Newton iteration, we use a metric induced by the Taylor expansion which can be defined in any Taylor category \cite{qualog23,TaylorMetric}.
	We define a family of operators on homsets $\M_n : \cat{C}(A, B) \to \cat{C}(A, B)$ mapping a morphism $f: A \to B$ to its \emph{$n$-th Taylor monomial} by evaluating the $n$-th derivative at $0$:
	\[
	\M_n(f)(a) :=  \frac{\mathsf{d}^nf(x)}{\mathsf{d} ^nx}(0)\cdot (a, \dots, a)
	\]
	This family of operators induces a pseudo-metric (the separation axiom does not necessarily hold) on each homset $d : \cat{C}(A, B) \times \cat{C}(A, B) \longrightarrow \cat{C}(A, B)$ given by
	\[
	(f,g) \longmapsto \begin{cases}
		2^{-k} & \text{where } k = \inf\{ n \in \NN \mid \M_n(f) \neq \M_n(g) \} \\
		0 & \text{ if for all } k \in \NN, \M_k(f) = \M_k(g)
	\end{cases} 
	\]
	similar to the Arnold-Nivat metric for $\lambda$-terms with truncations~\cite{arnold1980metric}.  In a Cartesian differential category, $d(f,g) = 2^{-k}$ also corresponds to the intuition $f(x) - g(x) = o(x^k)$ in analysis and we use this canonical distance for convergence analysis of optimization in our setting. The stronger triangle inequality holds for the Taylor pseudo-metric:
	\[
	d(f,h) \leq \max \{ d(f,g), d(g,h) \}.
	\]
	If $\cat{C}$ is a \emph{Taylor category} \cite{ehrhard2018introduction, kerjean2023taylor,ehrhard2023coherent,ehrhard2023coherentbimonad}, then the separation axiom holds and we obtain an ultra-metric \cite{qualog23,TaylorMetric}. 
	
	 For the convergence rate of Newton's method, we assume further that for $f: A \times X \to X$, $f(0,0) =0$ and that the second partial derivative of $f$ is nilpotent in its linear argument, meaning that there exists $p \geq 1$ such that for all $a,z,b$:
	 \[
	\left( \frac{\mathsf{d}f(a,v)}{\mathsf{d}v}(z)\right)^p\cdot b =0.
	 \]
	 similarly to \cite{pivoteau2008generation,pivoteau2012algorithms} where they require the Jacobian matrix to be nilpotent. For an additive map $h : X \to X$, if $h$ is nilpotent, then for all $k : Y \to X$, we have 
	 \begin{equation}
	 		\label{eqn:nilpotentrepetition}
	 	k \leq h^* \circ (k \ominus (h \circ k)).
	 \end{equation}
	To simplify the notation, we write again $y := \fix(f) :A \to X$, $z_n:= \Z_n(f): A \to X$, $z_{n+1} := \Z_{n+1}(f): A \to X$ and \[(\Diff_{z_n} [f](a))\cdot b := \frac{\mathsf{d}f(a,v)}{\mathsf{d}v}(z_n(a))\cdot b.\] Using (\ref{eqn:minusAdj}), we obtain that $y(a)\ominus z_{n+1}(a)$ is equal to
	\[
	(y(a) \ominus z_n(a)) \ominus \left(\left(\Diff_{z_n} [f](a)\right)^*\cdot (f(a,z_n(a)) \ominus z_n(a))\right)
	\]
	and by (\ref{eqn:nilpotentrepetition}), it is less than or equal to 
	\begin{gather*}
		(\Diff_{z_n} [f](a))^*\cdot \left((y(a) \ominus z_n(a)) \ominus (\Diff_{z_n} [f](a))\cdot (y(a) \ominus z_n(a)) \right)\\
		\ominus \left( (\Diff_{z_n} [f](a))^*\cdot (f(a,z_n(a)) \ominus z_n(a))\right)
	\end{gather*}
	which, by (\ref{eqn:minusAdd}) and (\ref{eqn:minusAdj}), is in turn is less than or equal to 
	\begin{gather*}
		 (\Diff_{z_n} [f](a))^*\cdot (((y(a) \ominus z_n(a)) \ominus (f(a,z_n(a)) \ominus z_n(a))) \\
		\ominus (\Diff_{z_n} [f](a))\cdot  (y(a) \ominus z_n(a)))
	\end{gather*}
	using (\ref{eqn:minusAdj}) and (\ref{eqn:minusTriang}), it is less than or equal to 
	\begin{gather*}
		(\Diff_{z_n} [f](a))^*\cdot (y(a) \ominus (f(a,z_n(a)) \\
	 +(\Diff_{z_n} [f](a))\cdot  (y(a) \ominus z_n(a)) ))
	\end{gather*}
Lastly, since $y(a)= f(a,y(a))$ and $z_n \leq y$, we can apply the Taylor formula (\ref{eqn:Taylorsubstract}) and obtain that:
\begin{gather*}
	y(a) = f(a,z_n(a)) + (\Diff_{z_n}[f](a)) \cdot (y(a)\ominus z_n(a)) \\
	+ \sum\limits_{k \geq 2} \frac{1}{k!} \frac{\mathsf{d}^kf(a,v)}{\mathsf{d}v^k}(z_n(a))\cdot (y(a)\ominus z_n(a))^k
\end{gather*}
Using (\ref{eqn:minusUnit}), we conclude that $y(a)\ominus z_{n+1}(a)$ is less than or equal to
\[
(\Diff_{z_n} [f](a))^*\cdot \left(	\sum\limits_{k \geq 2}\frac{1}{k!} \frac{\mathsf{d}^kf(a,v)}{\mathsf{d}v^k}(z_n(a))\cdot (y(a)\ominus z_n(a))^k\right)
\]

		We now state the main theorem of this section:
		\begin{theorem}
			\begin{enumerate}
				\item[]
				\item The family of Newton approximants $\{\Z_n\}_{n \in \omega}$ is an $\omega$-chain and its supremum $\Z := \bigvee_n \Z_n$ verifies $\Z = \Y$ and for all $n$, $\Z_n \leq \Y$, which means that the Newton approximants converge to the least fixpoint from below.
				\item 
				If we consider the canonical distance induced by the Taylor monomials, for $f : A \times X \to X$, assuming that $f(0,0)=0$ and that the second partial derivative of $f$ is nilpotent in its linear argument, then the convergence rate is quadratic\footnote{In the works on Newton's method for power series over $\omega$-continuous semi-rings (\eg \cite{kiefer2007convergence}), the terminology of exponential convergence is sometimes used as the number of accurately computed monomials of the power series doubles at each iteration. This is equivalent to quadratic convergence for the Taylor metric.}, for all $n$:
				\[
				d(\Z_{n+1}(f), \Y(f)) \leq(d(\Z_n(f), \Y(f)))^2.
				\]
			\end{enumerate}
		\end{theorem}

	\section{Conclusion and Future work}\label{sec:conclusion}
	
	In this paper, we provided a categorical framework for combining the theory of differentiation and the theory of fixpoints. We introduced the notion of Cartesian differential fixpoint categories (Section \ref{sec:fixdif}), and studied the case of Conway operators (Section \ref{sec:traced}), the relation between linearity and fixpoints (Section \ref{subsec:fix-linear}), and the closed setting (Section \ref{sec:closed}). We showed how many well-known examples are Cartesian differential fixpoint categories, such as weighted relations (Section \ref{ex:relations}), quantales profunctors (Section \ref{ex:quantale}), models induced by bifree algebras (Section \ref{sec:bifreefixdif}), and those induced by fixpoint objects (Section \ref{sec:lawvere}). We also showed how the Newton-Raphson scheme can be applied in our framework (Section \ref{sec:Newton}). We conclude with a brief discussion on potential future work and directions that build on the story of this paper. 
	
	Since the tangent bundle functor played such a crucial role, the first natural direction to take is to study tangent categories \cite{cockett2014differential} with fixpoint operators, where the defining axiom will be the tangent-fixpoint rule (\ref{eq:trace-fix}). Moreover, (\ref{eq:trace-fix}) should also be considered as an axiom for tangent categories with trace operators in the general monoidal setting. 
	
	There are also cofree Cartesian differential categories which are built via the Fa\`a di Bruno construction \cite{garner2020cartesian,cockett2011faa}. We aim to study parametrized fixpoint operators in these cofree models or understand how the Fa\`a di Bruno construction can be used to build new examples of Cartesian differential fixpoint categories. 
	
	Cartesian reverse differential categories \cite{cockett2019reverse} provide the categorical foundations of reverse differentiation, an important tool for automatic differentiation. Since fixpoints are also an important tool for automatic differentiation \cite{mazza2021automatic}, it would be worthwhile to understand the compatibility relation between parametrized fixpoint operators and reverse differential combinators. 
	
	Many optimization schemes are based on refinements of the  Newton-Raphson method. We would like to extend our construction to other iterative methods and study their application to approximate solutions of differential equations in our setting. 
	
	There are also many cases where we can only obtain fixpoints defined on a restricted domain which are used when we are interested in \emph{local} minima or maxima and not global solutions. We therefore aim to develop a theory of fixpoints for differential \emph{restriction} category \cite{cockett2011differential} which would allow us to capture local implicit function theorems that are more used in practice.

	\begin{acks}
		We thank the anonymous reviewers for their detailed feedback which helped improve the paper. We are also grateful to Masahito Hasegawa for discussions on Conway fixpoint operators. The first author is partially supported by the MUR FARE project "CAFFEINE" (Compositional and Effectful Program Distances). The second author is funded by an ARC DECRA (DE230100303) and an AFOSR Research Grant (FA9550-24-1-0008).
	\end{acks}
	
	\newpage

	\bibliographystyle{ACM-Reference-Format}
	\bibliography{biblio}


\end{document}